\documentclass[12pt]{amsart}    
\input{diagrams}    
\diagramstyle[scriptlabels]
\input{chords.sty}    

\newcommand{\nc}{\newcommand}

\nc{\FP}{\mathrm{P}}     
\nc{\FPM}{\mathbb{P}}    
\nc{\FPMA}{\mathbb{P}_+} 
\nc{\FMC}{\mathbb{M}}    
\nc{\FMCA}{\mathbb{M}_+} 
\nc{\Augm}{\mathbb{A}}   
\nc{\Tot}{\mathrm{Tot}}  
\nc{\Gr}{\mathrm{Gr}}    
\nc{\Const}{\mathrm{Const}} 
\nc{\PBW}{Poincar\'e-Birkhoff-Witt}   
\nc{\EE}{\mathfrak{U}} 
    
\def\Oone#1#2#3{\hbox{\normalfont1\hskip-#1em 
\special{ps:currentpoint currentpoint translate #2 #3 scale neg exch neg exch    
translate}I
}}    
\def\one{\Oone{0.27}{.85}{.96}}        
\nc{\Vect}{\mathtt{Vect}} 
\nc{\PROP}{PROP}    
\nc{\UA}{\mathbb{L}}     
\nc{\UL}{\mathbb{L}}    
\nc{\PL}{\Lie^C}    
\nc{\PLM}{\Lie^M }    
\nc{\Props}{\mathtt{Props} }    
\nc{\inn}{\mathrm{in}}    
\nc{\out}{\mathrm{out}}    
\nc{\CycOp}{\mathtt{CycOp}}    
\nc{\ModOp}{\mathtt{ModOp}}    
\nc{\AugModOp}{\mathtt{ModOp}_+}

\nc{\invtensor}{\underset{\leftarrow}{\otimes}}    
\nc{\rlarrows}{\begin{picture}(1,0.4)    
                \put(0,-0.1){\makebox(1,0.2){$\leftarrow$}}    
                \put(0,0.1){\makebox(1,0.2){$\to$}}    
              \end{picture}}    
\nc{\rra}{\begin{picture}(1,0.4)    
                 \put(0,-0.1){\makebox(1,0.2){$\lra$}}    
                 \put(0,0.1){\makebox(1,0.2){$\lra$}}    
              \end{picture}}    

\nc{\Left}{\mathbf L}  
\nc{\Right}{\mathbf R} 
\nc{\gr}{\operatorname{gr}}    
\nc{\Ho}{\operatorname{Ho}}    
\nc{\alt}{\operatorname{alt}}    
\nc{\Sym}{\operatorname{Sym}}    
\nc{\sym}{\operatorname{sym}}    
\nc{\id}{\operatorname{id}}    
\nc{\Der}{\operatorname{Der}}    
\nc{\im}{\operatorname{Im}}    
\nc{\Ker}{\operatorname{Ker}}    
\nc{\coker}{\operatorname{Coker}}    
\nc{\Col}{\operatorname{Col}}    
\nc{\ter}{\operatorname{ter}}    
\nc{\intl}{\operatorname{int}}    
\nc{\val}{\operatorname{val}}    
    
\nc{\TN}{{\cal N}}    
\nc{\Norm}{\operatorname{N}}    
\nc{\Nor}{\operatorname{N}}    
\nc{\Tor}{\operatorname{Tor}}    
\nc{\res}{\operatorname{res}}    
\nc{\Stab}{\operatorname{Stab}}    
\nc{\Hom}{\operatorname{Hom}}    
\nc{\chom}{\CH\!o\!m}    
\nc{\uhom}{\CH\!o\!m}    
\nc{\End}{\operatorname{End}}    
\nc{\holim}{\operatorname{holim}}    
\nc{\dirlim}{\underset{\rightarrow}{\lim}\,}    
\nc{\invlim}{\underset{\leftarrow}{\lim}\,}    
\nc{\com}{\operatorname{co}}    
\nc{\Th}{\operatorname{Th}}    
\nc{\Cech}{\check{C}}    
\nc{\Spec}{\operatorname{Spec}}    
\nc{\Spf}{\operatorname{Spf}}    
\nc{\MC}{\operatorname{MC}}    
\nc{\U}{\operatorname{U}}    
\nc{\Diff}{{\cal D}\mbox{\em iff}}    
\nc{\Mor} {{\cal M}or}    
\nc{\Ob}{\operatorname{Ob}}    
\nc{\cone}{\widehat}    
\nc{\Coder}{\operatorname{Coder}}    
\nc{\pr}{\operatorname{pr}}    
\nc{\diag}{\operatorname{diag}}    
    
\nc{\CHo}{{\cal{H}\mbox{\it{o}}}}    
\nc{\Mod}{{\mathtt{mod}}}           
\nc{\Modf}{{\mathtt{modf}}}           
\nc{\Modg}{{\mathtt{modg}}}           
\nc{\Ab}{{\mathtt {Ab}}}              
\nc{\Alg}{{\mathtt {Alg}}}     
\nc{\Hoalg}{{\mathtt {Hoalg}}}     
\nc{\Valg}{{\mathtt {Viral}}}     
\nc{\Algf}{{\mathtt {Algf}}}     
\nc{\Algg}{{\mathtt {Algg}}}     
\nc{\Coalg}{{\mathtt {Coalg}}}     
             
\nc{\dgc}{{\mathtt{dgc}}}    
\nc{\dgca}{{\mathtt{dgca}}}    
\nc{\dgcu}{{\mathtt{dgcu}}}    
\nc{\dgcuf}{{\mathtt{dgcuf}}}    
\nc{\dgcf}{{\mathtt{dgcf}}}    
\nc{\dgcg}{{\mathtt{dgcg}}}    
\nc{\dgcc}{{\mathtt{dgccc}}}    
    
\nc{\dgl}{{\mathtt{dglie}}}    
\nc{\dgla}{{\mathtt{dgla}}}    
\nc{\dglf}{{\mathtt{dglf}}}    
\nc{\dglg}{{\mathtt{dglg}}}    
\nc{\dga}{{\mathtt{dga}}}    
\nc{\art}{{\mathtt {art}}}    
\nc{\dgar}{{\mathtt {dgart}^{\leq 0}}}    
\nc{\simpl}{\Delta^{\op}\Ens}    
\nc{\Coll}{{\mathtt{Coll}}}    
\nc{\Kan}{{\mathtt {Kan}}}    
    
\nc{\Grp}{{\mathtt {Grp}}}    
\nc{\Cat}{{\mathtt {Cat}}}    
\nc{\Ens}{{\mathtt {Ens}}}    
\nc{\op}{{\operatorname{op}}}    
\nc{\Op}{{\mathtt{Op}}}    
    
\nc{\Lie}{{\mathtt{LIE}}}    
\nc{\Com}{{\mathtt{COM}}}    
\nc{\Ass}{{\mathtt{ASS}}}    
    
\nc{\pa}{\partial}    
    
\nc{\cal}{\mathcal}    
\nc{\CA}{\cal A}    
\nc{\CBB}{\cal B}    
\nc{\CB}{\cal B}    
\nc{\CC}{\cal C}    
\nc{\CDD}{\cal D}    
\nc{\CD}{\cal D}    
\nc{\CE}{\cal E}    
\nc{\CF}{\cal F}    
\nc{\CG}{\cal G}    
\nc{\CH}{\cal H}    
\nc{\CI}{\cal I}    
\nc{\CJ}{\cal J}    
\nc{\CK}{\cal K}    
\nc{\CL}{\cal L}    
\nc{\CM}{\cal M}    
\nc{\CN}{\cal N}    
\nc{\CO}{\cal O}    
\nc{\CP}{\cal P}    
\nc{\CQ}{\cal Q}    
\nc{\CR}{\cal R}    
\nc{\CS}{\cal S}    
\nc{\CT}{\cal T}    
\nc{\CU}{\cal U}    
\nc{\CV}{\cal V}    
\nc{\CW}{\cal W}    
\nc{\CZ}{\cal Z}

\renewcommand{\frak}{\mathfrak} 
\nc{\fa}{\frak a}    
\nc{\fg}{\frak g}    
\nc{\fk}{\frak k}    
\nc{\fh}{\frak h}    
\nc{\fm}{\frak m}    
\nc{\fn}{\frak n}    
\nc{\fA}{\frak A}    
\nc{\fC}{\frak C}    
\nc{\fI}{\frak I}    
\nc{\fS}{\frak S}

\nc{\nen}{\newenvironment}    
\nc{\ol}{\overline}    
\nc{\ul}{\underline}    
\nc{\lra}{\longrightarrow}    
\nc{\lla}{\longleftarrow}    
\nc{\Lra}{\Longrightarrow}    
\nc{\Lla}{\Longleftarrow}    
\nc{\Llra}{\Longleftrightarrow}    
\nc{\hra}{\hookrightarrow}    
\nc{\iso}{\overset{\sim}{\lra}}    
    
\nc{\notebox}[1]{{\begin{center}\fbox{\fbox{\parbox{11cm}{\small\sf     
                  #1}}} \end{center}}}     
\nc{\eproof}{ {\quad \ \hfill \fbox{\hspace{.2ex}} \\ \medskip}} 
\nc{\Thm}[1]{Theorem~\ref{#1}}    
\nc{\Prop}[1]{Proposition~\ref{#1}}    
\nc{\Lem}[1]{Lemma~\ref{#1}}    
\nc{\Cor}[1]{Corollary~\ref{#1}}    
\nc{\Conj}[1]{Conjecture~\ref{#1}}    
\nc{\Claim}[1]{Claim~\ref{#1}}    
\nc{\Defn}[1]{Definition~\ref{#1}}    
\nc{\Exa}[1]{Example~\ref{#1}}    
\nc{\Rem}[1]{Remark~\ref{#1}}    
\nc{\Note}[1]{Note~\ref{#1}}    
    
\nen{thm}[1]{\label{#1}{\bf Theorem.\ } \em}{}    
\nen{prop}[1]{\label{#1}{\bf Proposition.\ } \em}{}    
\nen{lem}[1]{\label{#1}{\bf Lemma.\ } \em}{}    
\nen{klem}[1]{\label{#1}{\bf Key Lemma.\ } \em}{}    
\nen{cor}[1]{\label{#1}{\bf Corollary.\ } \em}{}    
\nen{conj}[1]{\label{#1}{\bf Conjecture.\ } \em}{}    
\nen{claim}[1]{\label{#1}{\bf Claim.\ } \em}{}    
    
\nen{defn}[1]{\label{#1}{\bf Definition.\ } }{}    
\nen{exa}[1]{\label{#1}{\bf Example.\ } }{}    
    
\nen{rem}[1]{\label{#1}{\em Remark.\ } }{}    
\nen{note}[1]{\label{#1}{\em Note.\ } }{}    
\nen{exer}[1]{\label{#1}{\em Exercise.\ } }{}    
\nen{pf}{\begin{proof}}{\end{proof}}    
    
\begin{document}    
   
\title[Cyclic Operads and Chord Diagrams]    
{Cyclic Operads and Algebra of Chord Diagrams
}   
\author{Vladimir Hinich}    
\address{Department of Mathematics, University of Haifa,    
Mount Carmel, Haifa 31905,  Israel}    
\email{hinich@math.haifa.ac.il}    
\author{Arkady Vaintrob}    
\address{Institut des Hautes \'Etudes Scientifiques, 
 Bures-sur-Yvette, 91440 France; \quad  {\sf Address after September 1,
 2000: \ }  Department
of Mathematics, University of Oregon,  Eugene, OR 97403, USA}     
\email{vaintrob@math.uoregon.edu}    
\begin{abstract}    
We prove that the algebra $\CA$ of chord diagrams,   
the dual to the associated graded algebra of Vassiliev knot   
invariants,  is isomorphic to   
the universal enveloping algebra of a Casimir Lie algebra   
in a certain tensor  category (the \PROP\ for Casimir Lie algebras).   
This puts on a firm ground a known statement 
that the algebra  $\CA$  ``looks and behaves like a universal    
enveloping algebra''. An immediate corollary  of our result   
is the conjecture of~\cite{bgrt} on the Kirillov-Duflo isomorphism for 
algebras of chord diagrams.     
    
Our main tool is a general construction of a functor from the category   
$\CycOp$ of cyclic operads to the category $\ModOp$ of modular operads   
which is left adjoint to the ``tree part'' functor $\ModOp \to \CycOp$.     
The algebra of chord diagrams arises when this construction is   
applied to the  operad $\Lie$. Another example of this construction   
is Kontsevich's graph complex that corresponds to the operad   
$\Lie_\infty$ for homotopy Lie algebras.   
\end{abstract}    
\maketitle    
    
\section{Introduction}   
   
It is well known that the theory of knot invariants of finite type   
(or Vassiliev invariants) is closely connected with Lie algebras.   
The aim of this paper is to clarify this relationship   
and give it a precise formulation.   
   
Knot invariants of finite type are related to various areas   
of mathematics and theoretical physics and  have been in the focus of   
very intensive research ever since V.\,Vassiliev introduced them in 1989.    
One of the remarkable features of these invariants is that  they  can be
completely characterized in terms of combinatorial objects called 
weight systems. A weight system is a function on 
chord diagrams (configurations of    pairs of points on a   
circle) satisfying certain relations. The dual space $\CA$ of   
the space of weight systems is generated by chord diagrams and has   
a natural structure of a graded commutative and cocommutative Hopf algebra.   
It is called {\em the algebra of chord diagrams\/}.    
   
Bar-Natan~\cite{BN} and Kontsevich~\cite{Konts} discovered a construction   
that gives a family of Vassiliev invariants for every   
finite-dimensional Lie algebra with a metric (an invariant inner product).   
According to calculations of J.\,Kneissler~\cite{Kn},   
all invariants up to order $12$ come from Lie algebras. However    
for large orders this is not true, and there exists a more   
general construction that gives Vassiliev invariants which cannot   
be obtained from Lie algebras. This construction was found by the second   
author in~\cite{Vn} as a byproduct of an attempt to understand the   
relationship between Lie algebras and invariants of knots.       
It turned out that the theory of Vassiliev invariants naturally   
leads to the concept of a Yang-Baxter Lie algebra, an algebraic   
structure generalizing Lie algebras and Lie superalgebras, 
and they in turn can be used to produce knot invariants. 
Namely, every metric Yang-Baxter Lie algebra $\fg$   
gives an algebra homomorphism   
\begin{equation}\label{eq:atoug}   
W_\fg:\CA\to U(\fg) 
\end{equation}   
with values in the center of  the universal enveloping algebra of $\fg$.   
Every linear functional on $Z(U(\fg))$ produces a sequence of   
Vassiliev invariants.

The existence of homomorphism $W_\fg$ is not the only indication   
of the Lie-type behavior of the algebra $\CA$ of chord diagrams.   
In particular, $\CA$ can be described as the space generated   
by certain (ribbon) graphs (also known as Chinese characters)   
modulo some relations and the proof   
of this fact (see~\cite{BN}) is strikingly parallel to the proof of   
the \PBW\ theorem for Lie algebras.   This raises a natural question whether
the algebra $\CA$ can be  described as the center of the 
universal enveloping algebra of a 
Lie-type object which is universal with respect to morphisms~(\ref{eq:atoug}).   
   
In this paper we show that this is, indeed, the case   
and prove, in particular, that every Vassiliev invariant   
factors through the homomorphism $W_\fg$ for some $\fg$.   
   
This universal object, however, cannot be found among  
Yang-Baxter Lie algebras, and to define it we need to move to   
a slightly higher level of abstraction.   
\smallskip   
   
Metric Lie algebras can be defined not only in the category of   
vector spaces, but in arbitrary linear tensor category.   
One can construct a metric Lie algebra $\UL^M$   
in a certain category $\Lie^M$ universal in the sense that   
every metric Lie algebra in a tensor category $\CC$   
is the image of $\UL^M$ under a unique tensor functor $\Lie^M \to \CC$.   
   
A Lie algebra $\fg$ in a tensor category has a universal enveloping   
algebra $U(\fg)$   which is an ordinary associative algebra   
in the category of vector spaces.  The collection of   
maps~(\ref{eq:atoug}) can now be described as an algebra   
homomorphism $\CA\to U(\UL^M)$.   
This homomorphism however is not an isomorphism,   
and the starting point of our work was to understand to which extent    
it determines the algebra $\CA$.   
   
       An appropriate setup is provided  
by the more general notion of a Casimir Lie algebra   
(i.e.\ Lie algebra with an invariant symmetric two-tensor).   
Similarly to the case of metric Lie algebras,
we construct the universal Casimir Lie algebra   
$\UL^C$ in a certain tensor category $\Lie^C$ and a homomorphism $U(\UL^C) 
\to \CA$.        
One of the main results of the paper is that this map is an isomorphism.     
   
This result can be derived from injectivity of the natural map   
$U(\UL^C)\to U(\UL^M)$ which, in turn, follows from a detailed    
analysis of the categories $\Lie^C$ and $\Lie^M$.   
\smallskip   
   
To present the main constructions and results of the paper   
we need to use the language of operads and PROPs.   
   
The notion of operad appeared in algebraic topology in late 60s   
as a tool for describing algebraic operations on iterated loop spaces.     
In algebra, operads are used to encode classes of algebraic   
structures (algebras over operads) with collections of polylinear   
operations of type  $L^{\otimes n}\to L$   satisfying specific   
properties. In particular, there exists an operad called $\Lie$, such that    
$\Lie$-algebras in different tensor categories include ``usual''   
Lie algebras,  Lie superalgebras, as well as dg-Lie algebras and   
Yang-Baxter Lie algebras.    
For every operad $\CO$  there exists a free $\CO$-algebra with a   
given set of  generators.

An invariant inner product on an algebra $L$  over a 
field $k$ can be considered as   
an operation $b:L^{\otimes 2}\to k=L^{\otimes 0}$.   
The notions of \PROP\ and  algebras over \PROP s have been designed   
to handle this and more general operations of type   
$L^{\otimes n}\to L^{\otimes m}$. 
 
As in the case of operads, the formalism of PROPs allows to define     
algebras in arbitrary tensor categories.   
\PROP\ itself is a very small tensor category: its objects are natural   
numbers  $\text{\bf n}\in \mathbb{N}= \{\mathbf{0,1,2,} \ldots\}$    
with the tensor structure given by addition.  
   Mentioned above categories $\Lie^M$ and $\Lie^C$ are, in fact, PROPs, 
   such that the corresponding algebras are exactly metric and Casimir 
   Lie algebras.    
In particular, the objects   
$$   
 \UL^M=\mathbf{1} \in\Lie^M \quad \mathrm{and} \quad   
 \UL^C=\mathbf{1} \in\Lie^C   
$$   
are Lie algebras in the corresponding  categories and they   
can be viewed as the universal metric Lie algebra and the   
universal Casimir Lie algebra respectively.   
   
Similarly to the definition of an enveloping 
algebra of an algebra over an operad
(see~\cite{hla}, Sect.~3), one can 
consider two different versions of a universal  enveloping algebra of
the Lie algebra   $\UL^M \in\Lie^M$ (resp.\ $\UL^C \in\Lie^C$).   
The first, {\em internal\/} universal enveloping algebra,   
is an associative algebra in a certain   
extension of the category $\Lie^C$ (resp.\ $\Lie^M$). The second,   
{\em external\/}  universal enveloping algebra, is   
a genuine associative algebra. This algebra can be described as a collection 
of compatible endomorphisms of all representations of $\UL^M$ (resp.\ $\UL^C$).    
   
The homomorphisms~(\ref{eq:atoug}) giving    
Vassiliev invariants for arbitrary metric Yang-Baxter Lie algebras   
now can be interpreted as a single algebra homomorphism    
$$   
  W_{\UL^M}: \CA\to U^M   
$$   
from the algebra of chord diagrams to the external enveloping algebra   
of the universal metric Lie algebra $\UL^M \in\Lie^M$.   
 
A precise Lie-theoretical description  of the algebra of chord diagrams 
is obtained when we     replace metric Lie algebras by Casimir Lie algebras. 
The following theorem is the central result of the paper.   

\medskip   
   
\noindent
\begin{thm}{th:main}   
 There exists an algebra isomorphism   
$$ U^C          \iso \CA$$   
from the external enveloping algebra of the universal Casimir Lie algebra   
           $\UL^C \in\Lie^C$ to the algebra $\CA$  of chord diagrams.   
\end{thm}   
\medskip   
   
As an immediate consequence, we see that each Casimir Lie algebra   
$\fg$ gives rise to a homomorphism from $\CA$ to the center of    
the enveloping algebra of $\fg$ and that every Vassiliev invariant   
can be obtained from a Casimir Lie algebra in some tensor category.   
   
This explains the similarities between the algebra of chord diagrams   
and Lie algebras. In particular, the above-mentioned   
description of the algebra $\CA$ in terms of Chinese characters   
follows from the \PBW\ theorem for $U^C$.   
Another immediate corollary of this theorem is the  
conjecture of Bar-Natan, Garoufalidis, Rozansky, and Thurston~\cite{bgrt}   
on the existence of a Kirillov-Duflo-type  
isomorphism for algebras of chord diagrams.    
\medskip

A large part of our results and constructions for metric and Casimir   
Lie algebras remains true if the operad $\Lie$ is replaced by an   
arbitrary  {\em cyclic operad\/} $\CO$.  In particular,  we construct   
 a \PROP\ $\CO^C$ describing $\CO$-algebras endowed with an   
invariant symmetric two-tensor. It turns out that the operadic part    
$\FMC(\CO)$ of $\CO^C$ has an extra structure --- that of a   
{\em modular operad\/}.   
Moreover, the functor $\CO\mapsto\FMC(\CO)$ is a left adjoint to   
the natural ``tree part'' functor from modular operads to cyclic   
operads. We give an explicit construction of $\FMC(\CO)$ in terms   
of $\CO$.    
The \PROP\ $\CO^C$ can be expressed through $\FMC(\CO)$ in a    simple   
way. Similarly, the \PROP\ $\CO^M$ also can be described   
in terms of $\FMC(\CO)$.   
The explicit description of PROPs $\CO^C$ and $\CO^M$ allows one   
to deduce  the following result which is the key ingredient in our   
characterization  of the algebra $\CA$.   
\medskip

\noindent
\begin{thm}{}   
The natural morphism of \PROP s
$\CO^C\to\CO^M$ induces an isomorphism   
$$ \Hom_{\CO^M}(\text{\bf 0},\text{\bf 0})\otimes   
\Hom_{\CO^C}(\text{\bf 0},\text{\bf n}) \iso
\Hom_{\CO^M}(\text{\bf 0},\text{\bf n}).$$   
\end{thm}   
\medskip   
   
When $\CO=\Lie$ this theorem implies in particular that all   
Vassiliev invariants can be obtained from a metric Lie algebra.   
\bigskip   
   
Another interesting example is the case when $\CO=\Lie_\infty$, the    
operad for homotopy Lie algebras. In this case the space of morphisms   
in $\Lie_\infty^C$ coincides with Kontsevich's graph complex.   
\medskip   
   
The paper is organized as follows.   
The constructions and results  of the first 
part (Sections 2 -- 4) are valid for 
arbitrary cyclic operads. We believe that they may find other   
applications besides the ones we discuss in the second part of the paper   
(Sections 5 -- 7).   
   
Section~2 describes our conventions   
about tensor categories, operads, and PROPs. In Section~3 we study cyclic and 
modular operads and algebras over them. We construct various functors between 
the categories of cyclic and modular operads and PROPs. In Section~4 we 
prove some results on adjointness and isomorphisms for 
these functors.   In Section~5 we discuss various  versions of the notion of
universal enveloping  algebra for algebras over operads and PROPs.  We prove
that when $\CO=\Lie$\, internal enveloping algebras exist in a certain
extension of the tensor category $\CP$ and then 
study them in detail in the cases when $\CP=\Lie^M$ and $\CP=\Lie^C$. In
Section~6 we review  basic facts about Vassiliev knot invariants and their
relations  with Lie algebras. Finally, in Section~7, we present several
applications of the results of the previous sections. In particular, we show    
how results of Section~5 allow to describe the algebra $\CA$   
of chord diagrams as the external universal enveloping algebra   
of the universal Casimir Lie algebra $\UL^C$.                      
\medskip   
   
{\em Acknowledgments.\/}   
This work was started when both authors were visiting   
Max-Planck-Institut f\"ur Mathematik in Bonn and finished when   
we both were at MSRI in Berkeley. We express our gratitude to these   
Institutes for hospitality and financial support.    
   
Also we would like to thank     Silvio Levy for showing  
us how to  make $\one$.

\section{Preliminaries}   
   
\subsection{Tensor categories}   
\label{tensor}   
   
\   
   
By a {\em tensor category\/}   we understand a $k$-linear symmetric monoidal
category  (see~\cite{dm,del})  over a field $k$ of characteristic zero.  The
unit object in a tensor category will be usually  denoted by $\one$.   
For any object $A$ in a tensor category, 
the associativity constraint allows to define uniquely (up to a 
unique isomorphism) the tensor powers $A^{\otimes n}$ and the 
commutativity constraint gives a left action of the symmetric 
group $\Sigma_n$ on $A^{\otimes n}$.    
   
Recall the following definition.    
\subsubsection{}   
\begin{defn}{rigid}   
An object $A$ of a tensor   category $\CC$ is called {\em rigid}  
if there  exists an object $A^{\vee}\in\CC$ and a pair of morphisms   
$$  
\phi:\one\to A^{\vee}\otimes A,\quad \psi: A\otimes A^{\vee}\to\one~,
$$    
such that the compositions   
$$ 
 A^{\vee}\rTo{\phi\otimes\id}A^{\vee}\otimes A\otimes    
 A^{\vee}\rTo^{\id\otimes\psi}A^{\vee}
$$   
and   
$$ 
A\rTo^{\id\otimes\phi}A\otimes A^{\vee}\otimes A   
\rTo^{\psi\otimes\id}A
$$   
are the identities.   
   
The object $A^{\vee}$ is called {\em the dual of } $A$.    
The dual object, if it exists, is unique up to a unique isomorphism.   
The pair $(\phi,\psi)$ is called {\em an adjoint pair}. Given one   
of the morphisms $\phi$ or $\psi$, its adjoint, if it exists, is unique.   
\end{defn}   
   
\subsubsection{} \begin{exa}{eg1}   
Let $\CC$ be the category $\Vect$  of $k$-vector spaces. 
Then $V\in\CC$ is rigid if and only if\, $\dim V<\infty$.   
\end{exa}

\subsection{\PROP s and algebras over them}   
 
\ 
   
Here we recall some basic facts about \PROP s, operads, and   
algebras over them.   Standard  
        references for this material   are~\cite{bv,ca,may,ad}.   
See also~\cite{ke2}, 1.1, 1.2.   
   
\subsubsection{}   
\label{initprop}   
Denote by $\mathbb{S}$ the tensor category whose objects are non-negative   
integers  $\mathbf{0},\mathbf{1},\mathbf{2},\ldots$, and   
morphisms are given by   
\begin{equation}   
\Hom_{\mathbb S}(\mathbf{m},\mathbf{n})=\left\{   
   \begin{array}{ll}   
   \emptyset,& m\not=n \\   
   \Sigma_n, & m=n,   
   \end{array}   
\right.   
\end{equation}   
where  $\Sigma_n$ is the symmetric group on $n$ objects.   
Tensor product in $\mathbb{S}$ is given by the addition of numbers;   
commutativity constraint 
      $$ s_{mn}: \mathbf{m}\otimes\mathbf{n}\to \mathbf{n}\otimes\mathbf{m}$$ 
is defined by the shuffle    
$$s_{mn}\in\Sigma_{m+n},\quad s_{mn}(i)=\left\{\begin{array}{ll}   
i+n,& i\le m\\   
i-m,& i>m   
\end{array}\right..   
$$   
   
The category $\mathbb{S}$ is the simplest example of a \PROP\ 
(see~\cite{ca}).     
\subsubsection{}   
\begin{defn}{defprop-1}  A   \PROP\  is a tensor  category $\CP$  
with  $\Ob\CP=\{\mathbf{0},\mathbf{1},\mathbf{2},\ldots\}$   
and a tensor functor $\mathbb{S}\to\CP$ identical on objects and injective on    
morphisms.   
   
       For a \PROP\ $\CP$, we will write $\CP(m,n)$ instead of    
$\Hom_{\CP}(\mathbf{m},\mathbf{n})$.   
 
A morphism of \PROP s $f:\CP\to\CP'$ is a functor from $\CP$ to $\CP'$   
which is identical on $\mathbb{S}$.   
The category of \PROP s  will be denoted by $\Props$.   
\end{defn}

\subsubsection{}   
\begin{defn}{prop-alg}   
Let $\CP$ be a \PROP\ and let $\CC$ be a tensor category.    
A $\CP$-{\em algebra\/} in $\CC$ is a tensor functor $A:\CP\to\CC$.   
\end{defn}    
   
\subsubsection{}   
\begin{exa}{eg2}   
Let    $\Vect$ be the category of $k$-vector spaces.   
A $\CP$-algebra   in $\Vect$ is a vector space $V=A(\mathbf{1})$ together
with a compatible collection of operations 
$A(p): V^{\otimes m}\to V^{\otimes n}$,  for each   $p\in\CP(m,n)$.   
\end{exa}

\subsubsection{Operads}   
By an {\em operad\/} in this paper we understand  
a collection of vector spaces   $\CO=\{\CO(n)\},\ n\in\mathbb{N}$, 
           endowed   with a right action of the symmetric group 
$\Sigma_n$ on $\CO(n)$  and a collection of  composition maps   
 
\begin{equation}\label{def-comp-op}   
\gamma:\CO(n)\otimes\CO(m_1)\otimes\ldots\otimes\CO(m_n)\to\CO\left(\sum m_i\right)   
\end{equation}   
satisfying natural conditions of equivariance, associativity, and unity 
(see~\cite{may,hla}). 
 
An     {\em algebra\/} over an operad $\CO$ is a vector  
space $A$   with a collection of   operations  
$$ 
\CO(n)\otimes A^{\otimes n}\to A,  
$$ 
satisfying natural compatibility conditions.     
  
The category of operads will be denoted by $\Op$.   
   
\subsubsection{PROPs and operads} \label{seq:po} 
There exists a pair of adjoint functors      
$$ 
\#: \Props\to \Op  \text{ \ and \ } \FP: \Op \to \Props,   
$$   
where $\CP^{\#}(n)=\CP(n,1)$, and the left adjoint to $\#$ functor $\FP$ 
    is defined  by the formula   
\begin{equation}   
\FP(\CO)(m,n)=\bigoplus_f   \bigotimes_{i=1}^n \CO(|f^{-1}(i)|),   
\label{freeprop}   
\end{equation}   
where  $\CO$ is an operad and the direct sum is taken over all maps    
$$f:\{1,\ldots,m\}\to\{1,\ldots,n\}.$$

The notions of algebras over \PROP s and over operads are compatible: 
an algebra over an operad $\CO$ is the same as an algebra over the 
\PROP\ \, $\FP(\CO)$.    
   
By a map of an operad $\CO$ to a \PROP\ $\CP$ we understand   
a morphism of operads  $\CO\to\CP^{\#}$.

\section{Cyclic and modular operads}   
\label{c-m-operads}

\subsection{Metric and Casimir algebras}    
\label{m-c-algebras}    
    
\    
In this section we define two types of algebras over a cyclic operad.    
One type, that of metric algebras, is well known. The other one    
generalizes the concept of  
Lie algebra endowed with a  casimir  element.     
    
These two types of algebras are governed by two different     
\PROP s which will be the main objects of study in the paper.

        \subsubsection{Cyclic operads}   A {\em cyclic operad\/}  
        (see~\cite{gk}) is an   operad $\CO$    
with a right action of the symmetric group $\Sigma_{n+1}$    
on $\CO(n)$  extending the $\Sigma_n$-action  and satisfying the    
compatibility  condition~(\ref{CYC}) below.

Note that the symmetric group $\Sigma_{n+1}$ is generated  
     by the subgroup $\Sigma_n$    
(identified with the stabilizer of $0\in\{0,\ldots,n\}$) and by the cyclic    
permutation $\tau$ given by    $\tau(i)=i-1$ for $i>0$ and $\tau(0)=n$.

The operad structure on $\CO$ can be described in    
terms of composition operations  
$$\circ_i:\CO(m)\otimes\CO(n)\to\CO(m+n-1),\ i=1,\ldots,m,$$
corresponding to the insertion of an element of $\CO(n)$ as the $i$-th 
argument of an element of $\CO(m)$.    
\medskip

\begin{defn}{def:cycop} 
An operad $\CO$ with a collection of right $\Sigma_{n+1}$-actions on $\CO(n)$    
is called {\em cyclic\/} if    
\begin{equation}\label{CYC}   
(a\circ_1b)\tau=(b\tau)\circ_{n}(a\tau)\text{ for }a\in\CO(m),    
       b\in\CO(n).    
\end{equation}   
    
The category of cyclic operads will be denoted by $\CycOp$.    
\end{defn} 
 
\subsubsection{}    
\label{metric-alg}    
\begin{defn}{metric}    
Let $\CO$ be a cyclic operad. A {\em metric $\CO$-algebra\/}  
in a tensor category $\CC$  
is  an $\CO$-algebra $A\in \CC$ together with a  symmetric form     
$b:A\otimes A\to\one$ satisfying the following conditions.    
   
(i) \     
The form $b$ is non-degenerate, that is     
there exists a two-tensor $c:\one\to A\otimes A$ adjoint to $b$    
in the sense of Definition~\ref{rigid}.    
\medskip    
    
(ii) \ The form $b$ is $\CO$-invariant, that is    
the composition    
\begin{equation}    
\CO(n)\otimes A^{\otimes n+1}\to A\otimes A\overset{b}{\lra}\one    
\end{equation}    
is $\Sigma_{n+1}$-invariant.      
\end{defn}    

\medskip
 
In~\cite{gk} metric algebras are called cyclic algebras. 
 
\subsubsection{PROP for metric algebras}    
\label{metric-prop}    
    
The notion of a metric algebra gives rise to the following construction.    
    
For a cyclic operad $\CO$ let $\CO^M$ be the \PROP\
generated by the \PROP\ $\FP(\CO)$ (see~(\ref{freeprop}) ) 
and two elements $b\in\CO^M(2,0)$ and $c\in\CO^M(0,2)$ 
satisfying the following conditions.    
    
(i) The morphisms $b$ and $c$ are symmetric and mutually adjoint.    
    
(ii) ({\em invariance}) For each $f\in\CO(n)$ the composition    
\begin{equation}    
\label{O-inv-m}    
\mathbf{n}\overset{c\otimes\id}{\lra}\mathbf{2}\otimes\mathbf{n}=    
\mathbf{1}\otimes\mathbf{n} \otimes\mathbf{1}    
\overset{\id\otimes f\otimes\id}{\lra}\mathbf{1}\otimes\mathbf{1}\otimes    
\mathbf{1}\overset{\id\otimes b}{\lra}\mathbf{1}    
\end{equation}    
       is equal to $f\tau$.    
\smallskip    
    
The following simple result explains the meaning of $\CO^M$.    
    
\subsubsection{}    
\begin{lem}{alg(prop-m)}    
Metric $\CO$-algebras are precisely the algebras over     
the \PROP\ $\CO^M$.    
\end{lem}    
\eproof 
    
   For a cyclic operad $\CO$ there exists a canonical morphism    
\begin{equation}    
i_M^{\CO}:\CO\to\CO^M    
\label{i_M}    
\end{equation}    
 to the corresponding \PROP.    
This morphism gives the    
functor $${i_M^{\CO}}_*:\Alg(\CO^M)\to\Alg(\CO)$$     
that forgets the metric of a metric $\CO$-algebra.    
    
\subsubsection{}    
\label{casimir-alg}    
\begin{defn}{casimir}    
Let $\CO$ be a cyclic operad. A {\em Casimir $\CO$-algebra\/} in a 
tensor category  $\CC$   
is an $\CO$-algebra $A\in \CC$ together with a symmetric $\CO$-invariant  
two-tensor (called {\em casimir}) $c:\one\to A\otimes A$. 
 
       The condition of $\CO$-invariance     
means that the following composition  
\begin{equation} \label{eq:casinv} 
\CO(n)=\CO(n)\otimes\one^{\otimes n}\overset{c^{\otimes n}}{\lra}    
\CO(n)\otimes (A^{\otimes 2})^{\otimes n}\to A\otimes A^{\otimes n}=    
A^{\otimes n+1}    
\end{equation}    
is $\Sigma_{n+1}$-equivariant with respect to the 
standard $\Sigma_{n+1}$-action on    
$A^{\otimes n+1}$ given by 
$$ x\sigma=\sigma^{-1}(x)\text{ \ for \ }x\in A^{\otimes n+1},     
      \text{ and } \sigma\in\Sigma_{n+1}.$$    
\end{defn}     
    
\subsubsection{PROP for Casimir algebras}    
\label{casimir-prop}    
    
Similarly to~\ref{metric-prop} we construct a PROP 
responsible for Casimir algebras.  
 
   Let $\CO$ be a cyclic operad. Denote by $\CO^C$ the PROP 
generated by the \PROP\, $\FP(\CO)$ (see~\ref{seq:po}) and  
a symmetric element  
   $c\in\CO^C(0,2)$   satisfying the following invariance condition.    
    
 For each $f\in\CO(n)$ the diagram    
\begin{equation}\label{eq:invariance} 
\begin{diagram}    
\mathbf{n-1} & \rTo^{c\otimes\id} & \mathbf{2}\otimes\mathbf{n-1}=    
\mathbf{1}\otimes\mathbf{n} \\    
\dTo^{\id\otimes c} & & \dTo^{\id\otimes f}\\    
\mathbf{n-1}\otimes\mathbf{2}=\mathbf{n}\otimes\mathbf{1} &     
\rTo^{f\tau\otimes\id} & \mathbf{2}    
\end{diagram}    
\end{equation} 
is commutative.    
    
The following fact is an analog of    
Lemma~\ref{alg(prop-m)} for Casimir algebras. 
 
\subsubsection{}    
\begin{lem}{alg(prop-c)}    
Casimir $\CO$-algebras are precisely the algebras over the \PROP\ $\CO^C$.    
\end{lem}    
\eproof 
    
For    a cyclic operad $\CO$ there is a canonical morphism    
\begin{equation}    
i_C^{\CO}:\CO\to\CO^C    
\label{i_C}    
\end{equation}    
to the corresponding \PROP\ for Casimir $\CO$-algebras.    
This morphism  gives the functor    
$$  {i_C^{\CO}}_*:\Alg(\CO^C)\to\Alg(\CO)  $$    
that forgets the casimir of a Casimir $\CO$-algebra.    
\medskip 
 
It is easy to see that algebras with invertible casimirs are metric 
algebras.   
\subsubsection{}    
\begin{lem}{met=cas+fin}    
 Let $A$ be an algebra over a cyclic operad $\CO$ and 
 $$ 
b:A\otimes A\to\one \text{ \ and \ } 
c:\one\to A\otimes A 
$$  
be a pair of  symmetric mutually adjoint maps. 
Then $b$ satisfies the conditions (i), (ii) of~\Defn{metric} if and 
only if $c$ is $\CO$-invariant in the sense of~\Defn{casimir}.   
\end{lem}    
\eproof    
 
\ 
 
This lemma gives a functor 
$$ 
\Alg(\CO^M)\to\Alg(\CO^C) 
$$ 
that commutes with $i^{\CO}_{M*}$ and  $i^{\CO}_{C*}$ and is induced by 
a morphism of \PROP s 
\begin{equation}\label{eq:icm} 
\CO^C\to \CO^M 
\end{equation} 
commuting with $i^{\CO}_M$ and  $i^{\CO}_C$. 
    
The goal of this section is to give           a detailed description of the
relationship between the \PROP s  $\CO^M$ and $\CO^C$.

\subsection{Coordinate-free language} \label{sec:coord-free}
\    
    
In this paper, when dealing with tensor categories and operads we will use a    
``coordinate-free'' language  of~\cite{dm}. It allows one to     
     hide some ``ugly'' part of structure (the associativity and     
commutativity constraints, action of  symmetric group, etc.) inside 
the category of finite sets. In this subsection we recall the basic     
definitions and reformulate the notion of a cyclic operad in the new    
language. In the following  subsection we will use this language to describe
modular operads.     
    
 The following is a coordinate-free definition of tensor category     
  (see \cite{dm}, Prop. 1.5).       
\subsubsection{}    
\begin{defn}{def:tns}  
A      {\em tensor category\/}  $\CC$ is a category with  
functors 
$$ 
\bigotimes_I:\CC^I\to\CC   \ : \ (X_i,i\in I) \mapsto \bigotimes_{i\in I}X_i 
$$    
and functorial isomorphisms    
$$ 
\chi(\alpha): \bigotimes_{i\in I} X_i \iso \bigotimes_{j\in J}    
\left(\bigotimes_{i\in\alpha^{-1}(j)} X_i \right) 
$$    
defined for each finite set $I$ and  each map $\alpha: I\to J$ of finite sets.  
The functors $\displaystyle\bigotimes_I$ and isomorphisms $\chi(\alpha)$ 
have to satisfy the following conditions. 
 
(i) If $I$ consists of a single element, then $\otimes_{i\in I}$ is    
the identity functor, and for any map $\alpha$ between  
one-element sets
$\chi(\alpha)$ is the identity automorphism of the identity functor.    
    
(ii) The functors $\chi(\alpha)$ satisfy a natural associativity condition 
relating $\chi(\beta\circ\alpha)$ with $\chi(\beta)$ and $\chi(\alpha)$  
for any pair of maps $I\rTo^\alpha J \rTo^\beta K$  of finite sets.       
\end{defn}

\subsubsection{\PROP s in the coordinate-free language} 
In the new language the initial \PROP\ $\mathbb{S}$ (see~\ref{initprop})    
is replaced by the groupoid of finite sets with tensor product given by the    
operaton of disjoint union. 
As a tensor category it is equivalent to the category of $\mathbb{S}$ 
of~\ref{initprop} and  
we will denote it by the same symbol. 
 
A \PROP\ in this setting is defined as a tensor  
category $\CP$  
with a tensor functor from $\mathbb{S}$  to $\CP$ bijective on objects 
and injective on morphisms.    
    
\subsubsection{Cyclic operads in the coordinate-free language}    
\label{ss:cyclic-operad}    
         Denote by  $\mathbb{S}^*$ the groupoid of {\em non-empty\/}  
  finite sets.     
A      {\em cyclic operad\/}  
is a functor    
$$    
 \CO: \mathbb{S}^* \to \Vect
$$    
with a collection of functorial composition operations 
\begin{equation}\label{eq:comp}    
\circ_{xy}: \CO(X) \otimes \CO(Y) \to \CO(X\sqcup Y \setminus \{x,y\})     
\end{equation}    
defined for each pair  
of pointed sets $(X,x)$ and $(Y,y)$ 
satisfying the following conditions.    
    
\begin{enumerate} 
\item[(i)] 
 ({\em commutativity}) Operations $\circ_{xy}$ and $\circ_{yx}$ coincide     
after  canonical identification of $X\sqcup Y$ with $Y\sqcup X$.    
    
\item[(ii)] ({\em associativity}) For $x\in X,\quad y,y'\in Y,\ y\not=y',    
\quad z\in Z$, the following diagram is commutative    
    
$$    
\begin{diagram}    
\CO(X)\otimes \CO(Y)\otimes \CO(Z) & \rTo^{\circ_{xy}\otimes\id_{\CO(Z)}} &   
\CO(X\sqcup  Y\setminus\{x,y\}) \otimes \CO(Z) \\    
\dTo_{\id_{\CO(X)}\otimes\circ_{y'z}} & & \dTo_{\circ_{y'z}} \\    
\CO(X)\otimes \CO(Y\sqcup Z\setminus\{y',z\}) & \rTo^{\circ_{xy}}    
& \CO(X\sqcup Y\sqcup Z\setminus\{x,y,y',z\}) .    
\end{diagram}    
$$    

\medskip
    
\item[(iii)] ({\em unity}) For every two-element  
      set $\{p,q\}$,     
there is a distinguished element $I_{pq} \in \CO(\{p,q\})$, such that     
for any $x\in X$ and $ a\in \CO(X)$, the elements    
$a$ and $a\circ_{xy}I_{yz}$ coincide after identifying $X$ with $X\cup\{z\}    
\setminus\{x\}$.    
\end{enumerate} 
    
\medskip

The space $\CO(X)$ should be viewed as a set of ``relations'' with arguments     
labelled by $X$.  
The composition $\circ_{xy}$ corresponds to the 
operation of gluing the sets $X$ and $Y$ along the points $x$ and $y$.    
    
\subsection{Modular operads}    
\ 
 
Roughly speaking, modular operads are cyclic operads where gluing 
operations  similar to~(\ref{eq:comp})  are allowed for arbitrary 
non-empty subsets $U\subset X$  and $V\subset Y$. We will use two 
different kinds of modular operads (see 
definitions~\ref{ss:mod-operad} and~\ref{amod-op}).

Modular operads were introduced by Getzler and Kapranov in~\cite{gkm}.     
Our definitions differ slightly from the one     
given in~\cite{gkm}  --- see~\ref{app.mod.op} for a comparison.

\subsubsection{}    
\begin{defn}{ss:mod-operad}    
A {\em modular operad\/} is a collection of functors    
$$\CM^n:\mathbb{S}^* \to \Vect,\quad n=0,1,\ldots,$$
with composition  operations    
\begin{equation}\label{eq:mod_comps}    
\circ_f: \CM^m(X) \otimes \CM^n(Y) \to \CM^{m+n+d-1}(X\sqcup Y    
\setminus(U\sqcup V)) 
\end{equation}    
defined for each bijection  
$U\rTo^f V$ between non-empty $d$-element subsets  
$U\subset X$ and $V\subset Y$ 
with 
\begin{equation} \label{eq:nonempty}    
X\sqcup Y \setminus(U\sqcup V)\not=\emptyset    
\end{equation}    
satisfying the following conditions.    
\smallskip    
    
\begin{enumerate} 
\item[(i)] 
({\em commutativity}) \ Operations $\circ_f$ and $\circ_{f^{-1}}$ 
coincide after the canonical identification of $X\sqcup Y$ with $Y\sqcup X$.     
    
\item[(ii)] 
({\em associativity}) \ For $i=1,2,3$, let $X_i$ be a non-empty finite set    
with two {\em disjoint\/ } subsets $U_{ij}\subseteq X_i, \
j\in\{1,2,3\}\setminus \{i\}$.  
Let, in addition, $f_{ij}:U_{ij}\to U_{ji}$     
be bijections satisfying  $f_{ij}=f^{-1}_{ji}$.  
 
If $U_{13}=U_{31}=\emptyset$ and  
$U_{12}$ and $U_{23}$ are non-empty, then    
\begin{equation}    
\label{eq:ass1}    
\circ_{f_{23}}(\circ_{f_{12}}\otimes\id_{X_3})    
=    
\circ_{f_{21}}(\circ_{f_{23}}\otimes\id_{X_1}).    
\end{equation}    
    
If all the subsets $U_{ij}$ are     
non-empty, then the following three     
maps from $\displaystyle\bigotimes_{i=1}^3\CM(X_i)$  
to $\CM\Bigl(\bigsqcup\limits_i (X_i\setminus\bigcup\limits_j U_{ij})\Bigr)$    
coincide:    
\begin{equation}\label{eq:ass2}    
\circ_{f_{13}\sqcup f_{23}}(\circ_{f_{12}}\otimes\id_{X_3})=    
\circ_{f_{12}\sqcup f_{32}}(\circ_{f_{13}}\otimes\id_{X_2})    
=\circ_{f_{21}\sqcup f_{31}}(\circ_{f_{23}}\otimes\id_{X_1}).    
\end{equation}    
    
\item[(iii)] ({\em unity}) \ See~\ref{ss:cyclic-operad}(iii).
\end{enumerate} 
    
\end{defn}    
 
\smallskip 
    
The category of modular operads  
will be denoted by $\ModOp$.    
     
\medskip 
 
The following proposition gives an equivalent definition 
of a modular operad (see~\cite{gkm}, 3.4 -- 3.7).    
    
\subsubsection{} 
\begin{prop}{mod-operad2} A modular operad is a graded 
cyclic operad  
$$ 
   \CM=\bigoplus_{n\ge 0} \CM^n:\mathbb{S}^*\to\Vect 
$$ 
 endowed  with {\em     contraction\/}  operations    
$$    
c_{xy}:\CM^n(X)\to \CM^{n+1}(X\setminus \{x,y\}),    
\quad x,y\in X, x\ne y, X\not=\{x,y\}    
$$    
satisfying the following properties. 
    
\begin{enumerate} 
\item[(i)] 
 $c_{xy}=c_{yx}$.    
 
\item[(ii)]  
If $x,y,z,t\in X$ are four distinct elements, then  
the contractions   $c_{xy}$ and $c_{zt}$ commute.    
    
\item[(iii)] 
Let $x_1\ne x_2\in X,\quad y_1\ne y_2\in Y$ 
and $X\sqcup Y\ne\{x_1,x_2,y_1,y_2\}.$  
Then  the operations $ c_{x_1,y_1}\circ_{x_2,y_2}$ and 
$c_{x_2,y_2}\circ_{x_1,y_1}$ from 
$ \CM(X)\otimes \CM(Y)$ to $\CM(X\sqcup Y\setminus\{x_1,x_2,y_1,y_2\})$ 
coincide.     
\end{enumerate} 
 
\end{prop}    
    
\begin{pf}    
Suppose $\CM=\{\CM^n:\mathbb{S}^*\to\Vect\}$ is a modular operad in 
the sense  of~\ref{ss:mod-operad}. Then one can define the contraction 
operation     
$$ c_{xy}:\CM^n(X)\to \CM^{n+1}(X\setminus\{x,y\})$$    
as the composition with the identity $I_{x'y'}\in \CM^0(\{x',y'\})$     
under the map $f$ sending $x'$ to $x$ and $y'$ to $y$.     
 Property~(i) follows from \ref{ss:mod-operad}(i), property~(ii) 
 follows from \ref{ss:mod-operad}(ii) with $X_2=X$, $X_1=\{x,y\}$, 
 and $X_3=\{z,t\}$.   Property~(iii)  follows from the fact that both
 compositions   coincide with $\circ_f$,  
where  $f: \{x_1,x_2\} \rTo \{y_1,y_2\}: f(x_i)=y_i$. 
 
Consider now a graded cyclic operad $\CM$ endowed with a collection of    
contractions satisfying properties (i)--(iii) above.  
    Define  compositions 
$$  
\circ_f: \CM^m(X) \otimes \CM^n(Y) \to   
\CM^{m+n+d-1}(X\sqcup Y \setminus(U\sqcup V)) 
$$     
as follows. Choose $u\in U$ and define $\circ_f$ to be the composition    
of $\circ_{u,f(u)}$ with the contractions $c_{v,f(v)}$  
for all $v \in 
U\setminus\{u\}$.  
The result does not depend on the choice of $u\in U$ and on the order 
of the contractions by the properties~(ii) and~(iii) of contractions.     
\end{pf}    
 
\medskip 
    
       The following  
is a non-graded version of the notion of a modular operad.         
\subsubsection{}    
\label{ngr}    
\begin{defn}{}    
A {\em non-graded}     
modular operad is a cyclic operad $\CM$ together with contraction 
operations     
$$    
c_{xy}:\CM(X)\to \CM(X\setminus \{x,y\}),\quad 
                   x,y\in X,\ x\ne y,\ X\not=\{x,y\}     
$$    
satisfying properties (i)--(iii) of Proposition~\ref{mod-operad2}
\end{defn}   
    
The category of non-graded modular operads will be denoted  
by $\ModOp^{ngr}$. 
    
\subsection{Standard functors}    
\label{standard-f}    
\ 
 
Here we will construct several functors connecting 
various categories of operads and \PROP s.

\subsubsection{From modular operads to \PROP s}  \label{seq:newprop} 
We start with a natural construction that associates a  \PROP\ to a 
modular operad.    
\medskip

\begin{prop}{pr:propAP}    
There exists a functor  
$$  
\FPM:\ModOp\to\Props 
$$ 
with $\FPM(\CM)(X,Y)$  
given for $X,Y\in\mathbb{S}$ by the formula 
\begin{equation}\label{eq:fpm}    
\FPM(\CM)(X,Y) = \bigoplus_{\substack{X = \coprod_{i\in I}X_i \\     
   Y=\coprod_{i\in I}Y_i\\ Y_i\ne \emptyset }}    
   \bigotimes_{i\in I}\CM(X_i\sqcup Y_i).    
\end{equation}    
\medskip    
\end{prop}    
\begin{proof}    
It is sufficient to define a composition for elements     
$f\in\FPM(\CM)(X,Y)$ of the special type $f=\alpha\sqcup I$, where     
$X=X'\sqcup X''$, \ $Y=Y'\sqcup Y''$, \  $I=\sqcup_{x\in X''} I_{x,\phi(x)}$,
\ $\phi:X''\to Y''$ is a bijection, and $\alpha \in\CM(X'\sqcup Y')$.    
 
     Composition of elements of this type is defined using the composition
in the modular operad $\CM$.     
\end{proof}

\subsubsection{From modular to cyclic operads}    
\label{zero-component}    
If $\CM$ is a modular operad, its zeroth component $\CM^0$ is a cyclic    
operad. This gives a functor    
\begin{equation}\label{eq:grzero} 
\Gr^0:\ModOp \to \CycOp.    
\end{equation}

\medskip 
 
For non-graded modular operads there is a natural forgetful functor     
\begin{equation}\label{eq:forget1}    
          \#_{ngr}: \ModOp^{ngr}\to\CycOp: \quad A \mapsto A^\#.  
\end{equation}

  \subsubsection{From non-graded modular operads to graded and back}    
        There exists a pair of adjoint functors    
$$\Tot:\ModOp\to \ModOp^{ngr} 
$$ 
and 
$$ \Const: \ModOp^{ngr} \to \ModOp  
$$   
defined by the following formulas 
\begin{equation} 
\Tot(\{\CM^n\})=\bigoplus_n\CM^n \mbox{ \ and \ } \Const(\CM)^n=\CM.    
\end{equation} 
 
\medskip 
 
The composition $\Gr^0 \circ \Const$ is     
isomorphic to the forgetful functor $\#_{ngr}$~(\ref{eq:forget1}). 
    
\subsection{Augmented cyclic and modular operads}    
\label{augm}    
\    
 
We will    need the       
following variation of the notions of cyclic and modular operads. 
 
\subsubsection{}    
\begin{defn}{acyc-op}    
An      {\em augmented cyclic operad\/}  
is a functor on the groupoid of finite sets 
$$    
 \CO: \mathbb{S} \to \Vect, \quad     
$$    
         endowed with a collection of compositions 
\begin{equation} 
\circ_{xy}: \CO(X) \otimes \CO(Y) \to \CO(X\sqcup Y \setminus \{x,y\})     
\end{equation}    
defined for each pair $x\in X, \ y\in Y$ 
satisfying the commutativity, associativity, and unity conditions 
(i), (ii), and (iii) of~\ref{ss:cyclic-operad}. 
\end{defn} 
 
\subsubsection{}    
\begin{rem}{} To define a structure of     
augmented cyclic operad on a cyclic operad $\CO$, one has to choose    
a graded vector space $\CO(\emptyset)$  
and to define an operation 
$$ 
\circ_{xy}:\CO(\{x\})\otimes \CO(\{y\}) \rTo \CO(\emptyset). 
$$ 
    
In particular, every cyclic operad can be considered     
as an augmented cyclic operad with $\CO(\emptyset)=0$.    
\end{rem}

\subsubsection{}    
\begin{defn}{amod-op}    
An {\em augmented modular operad\/} is   
 a collection of functors     
$  \CM^n:\mathbb{S} \to \Vect,\ n=0,1,\ldots, $  
with a collection of compositions    
$$    
\circ_f: \CM^m(X) \otimes \CM^n(Y) \to \CM^{m+n+d-1}(X\sqcup Y 
\setminus(U\sqcup V))  
$$     
defined for each bijection  
$U\rTo^f V$ between non-empty $d$-element subsets  
$U\subset X$ and $V\subset Y$ satisfying the conditions 
(i)--(iii) of~\ref{ss:mod-operad}.     
\end{defn}    
\medskip    
       
The category of augmented operads will be denoted $\AugModOp$.    
\medskip    
    
The following is a 
    version of~\Prop{mod-operad2} for augmented modular operads.    
 
\subsubsection{}    
\begin{prop}{amod-op2}    
An augmented modular operad can be defined as 
an augmented graded cyclic operad  
  $$\CM^n:\mathbb{S}\to\Vect,\quad n=0,1,\ldots,$$    
endowed with contraction operations 
$$    
c_{xy}:\CM^n(X)\to \CM^{n+1}(X\setminus \{x,y\}),    
\quad x,y\in X, \ x\ne y,     
$$    
satisfying properties (i)--(iii) of~\Prop{mod-operad2},    
        where in property (iii) we do not require that 
 $X\sqcup Y\ne\{x_1,x_2,y_1,y_2\}$.   
\end{prop}    
               \eproof

\subsubsection{\PROP\ from an augmented modular operad}  
 \label{seq:a-newprop}     
Similarly to~\ref{seq:newprop}, we consider 
a functor $\FPMA$ that associates a \PROP\ to 
an augmented modular operad.     
 
Let $\CM$ be an augmented modular operad.     
Define a collection of vector spaces  
$\FPMA(\CM)(X,Y)$ for $ X,Y\in\mathbb{S}$ by the  formula    
\begin{equation} 
\FPMA(\CM)(X,Y) = S(\CM(\emptyset))\otimes    
\Bigl(\bigoplus_{\substack{X = \coprod_{i\in I}\!X_i \\     
    Y=\coprod_{i\in I}Y_i\\ X_i\sqcup\!Y_i\ne\emptyset 
    }}\bigotimes_{i\in I} \CM(X_i\sqcup Y_i) \Bigr), 
\end{equation} 
where  
$S(V)$ is the symmetric algebra of the vector space $V$. 
\medskip    
    
An argument similar to the proof of Proposition~\ref{seq:newprop} 
shows that the assignment $\CM\mapsto\FPMA(\CM)$ gives 
a functor    
$$ \FPMA:\AugModOp\to\Props.$$  
 
\medskip 
 
We will also need the following forgetful functor 
\begin{equation}\label{eq:forget23}    
              \#: \AugModOp\to\ModOp : \quad \CM \mapsto \CM^\#    
\end{equation}    
given by the formula 
$$  
\CM^\#(X)=\CM(X), \ \mathrm{for} \ X\ne \emptyset. 
$$ 
    
\subsubsection{Comparison with the definition of~\cite{gkm}}    
\label{app.mod.op}    
    
The modular operads of Getzler-Kapranov~\cite{gkm}    
   are in our terminology augmented modular operads     
$\CM$ satisfying the following additional {\em stability}  
requirements:    
    
\begin{itemize}    
\item $\CM^0(\emptyset)=\CM^1(\emptyset)=0;$    
\item $\CM^0(\{x\})=0;$     
\item $\CM^0(\{x,y\})=k\cdot I_{xy}.$    
\end{itemize}

\subsection{Results}    
    
\    
    
In this section we formulate the main results of the first part of the    
paper. They claim the existence of some adjoint functors to the    
standard functors defined in~\ref{standard-f} and \ref{augm}     
together with various relations between these functors.  
Proofs of these results
will be given in the next section.      
    
The following diagram shows  
relevant categories and functors.    
 The triangles formed by solid arrows  
will be commutative. 
    
\begin{equation}\label{eq:funcdiag} 
\begin{diagram}   
  & &\CycOp & & & & \\   
  & \ldTo^{\CO\mapsto \CO^C} & \uDotsto^{\Gr^0}\dTo_{\FMC}    
     &\rdTo(4,2)^{\CO\mapsto \CO^M} & & & \\   
 \Props&\lTo^{\FPM}&\ModOp&\rTo^{\Augm}&\AugModOp & \rTo^{\FPM_+} & \Props \\    
  & &\uDotsto^{\Const}\dTo_{\Tot} & & & & \\    
  & &\ModOp^{ngr} & & & &     
\end{diagram}    
\end{equation}

\subsubsection{}   
\begin{thm}{t:left-adj}    
The zero-component functor~(\ref{eq:grzero}) $$\mathrm{Gr}^0:\ModOp\to\CycOp$$
admits a left adjoint functor    
\begin{equation} \label{eq:funm} 
 \FMC: \CycOp\to\ModOp.  
\end{equation} 
\end{thm}    
    
\subsubsection{}    
\begin{cor}{cor:left-adj}    
The functor  $$ \Tot\circ\FMC: \CycOp\to\ModOp^{ngr}$$     
is a left adjoint to the forgetful functor~(\ref{eq:forget1}) 
          $$\#_{ngr}:\ModOp^{ngr}\to\CycOp.$$ 
\end{cor}    
 
\subsubsection{}    
\begin{rem}{rem:freemod}    
Let $V:\mathbb{S}^*\to \Vect$ be a functor and $F(V)$ be the free    
cyclic operad generated by $V$ (see~\cite{gk}).     
Since by \Thm{t:left-adj} the functor $\FMC$ is left adjoint to    
the zero-component functor $\Gr_0$ and the functor $F:V\mapsto F(V)$    
is left adjoint to the forgetful functor $\#_{ngr}$,  
the functor 
$\FMC\circ F$ is left adjoint to the composition of the 
functor  $\#_{ngr}$ 
with $\Gr^0$. 
This shows  
that $\FMC(F(V))$  
can be considered as the free modular operad    
generated by $V$~(see \cite{gkm}).    
\end{rem}    
    
\subsubsection{}    
\begin{thm}{coc=}    
The functors $\FPM \circ \FMC$ and $\CO\mapsto \CO^C$ 
from $\CycOp$ to $\Props$ 
are isomorphic, i.e.\ the left solid triangle of the 
diagram~(\ref{eq:funcdiag}) is commutative. 
\end{thm}    
    
\subsubsection{}    
\begin{thm}{augmentation}    
The forgetful functor    
$$\#:\AugModOp\to\ModOp$$    
admits a left adjoint  
\begin{equation}\label{eq:funca} 
\Augm:\ModOp\to\AugModOp~. 
\end{equation} 
\end{thm}

\subsubsection{}    
\begin{thm}{com=}    
Let $\FMCA$ denote the composition of functors 
$\Augm\circ\FMC$. 
 
The functors $\FPMA \circ \FMCA$ and $\CO \mapsto \CO^M$    
from $\CycOp$ to $\Props$ 
are isomorphic.  
In other words, the right solid triangle of the 
diagram~(\ref{eq:funcdiag} is commutative.  
\end{thm}    
    
\subsubsection{}    
\begin{cor}{coc-vs-com}    
Let $\CO$ be a cyclic operad. For each $X\in\mathbb{S}^*$  
there exists a  natural isomorphism of vector spaces     
\begin{equation}    
  \label{eq:com-vs-coc}    
 \CO^M(\emptyset,\emptyset)\otimes\CO^C(\emptyset,X)    
 \iso\CO^M(\emptyset,X). 
\end{equation}    
\end{cor}    
    
\subsubsection{}    
\begin{cor}{cor:inj}    
Let $\CO$ be a cyclic operad. 
The natural map    
$$\CO^C(0,n)\to\CO^M(0,n)$$    
is injective for all $n\ge 0$. 
\end{cor}    

\section{Proofs}   
\label{sec:proof}   
   
In this section we prove 
the results formulated in Sections~\ref{t:left-adj}--\ref{cor:inj}.     
The technical heart of the proof is the fact that 
the operadic part of the \PROP\ $\CO^C$ 
 admits a natural structure  of a modular operad.   
This is established in~\ref{mod-from-cyc}. 
 
Then in~\ref{therest} we show how to deduce the statements 
of \ref{t:left-adj} ---~\ref{cor:inj} from this fact.

\subsection{A modular operad from a cyclic operad}   
\label{mod-from-cyc}   
 
\ 
 
Let $\CO$ be a cyclic operad.   
Consider  a family of  vector spaces  
labelled 
by pairs $(X,x),$  
\ $X\in\mathbb{S}^*,\quad x\in X$  given by 
\begin{equation}\label{eq:mxx} 
  \CM_x(X)=\CO^C(X\setminus\{ x\}, \{x\}),   
\end{equation} 
                          where 
the right-hand side is understood in the PROP sense.   
 
This is just the                    {\em operadic part\/} 
 of the \PROP\ $\CO^C$.   
\medskip   
   
Our goal is to introduce a structure of a modular operad on this collection   
of spaces $\CM$. In particular we will canonically identify   
$\CM_x(X)$ for different $x\in X$. The resulting modular operad will 
be denoted $\FMC(\CO)$.   
   
\subsubsection{Grading}  
Recall that $\CO^C$ is defined in~\ref{casimir-prop} 
as the \PROP\  generated over $\FP(\CO)$ by the casimir $c\in\CO^C$, 
subject to     relations~(\ref{eq:invariance}). 
These relations  are homogeneous with respect to   the number of casimirs,  
therefore,  the space $\CM_x(X)$ obtains a natural grading 
$$  
\CM_x(X)=\bigoplus_{n\ge 0}\CM^n_x(X), 
$$  where   
$\CM^n_x(X)$ is the space generated by the compositions  
\begin{equation}  \label{cnf}   
               f \circ c^{\otimes n}, 
\end{equation} 
where 
$c^{\otimes n}:\emptyset\to Y\sqcup Y', \text{ with } |Y|=|Y'|=n $ 
        and 
$  f\in\CO(X\sqcup Y\sqcup Y'). $   
   
Note that  
since the presentation~(\ref{cnf}) is not unique, we cannot   
use it to identify    
  $\CM_x(X)$ for different $x$.   
The degree zero part of $\CM_x(X)$ coincides with $\CO(X)$.   
 
We will introduce the structure of a modular operad on the collection   
of spaces $\CM^n_x(X)$  by induction on degree.   
\medskip   
   
For a pointed set $(X,x)$ denote by $\widehat{X}$ 
the set $X\sqcup\{y,y'\}$ and define a morphism 
\begin{equation}   
\label{c:def}   
 c=c_{yy'}: \CM_x^n(\widehat{X}) \to \CM^{n+1}_x(X)    
\end{equation}   
as the composition with the casimir $c_{yy'}\in \CO^C(\emptyset,\{y,y'\})$   
in the \PROP\ $\CO^C$.   
\medskip   
   
\subsubsection{}   
 
\begin{lem}{properties-of-c}   
 
\ 
 
1.     The map $c_{yy'}$ is a surjection for $n\ge 0$.   
   
2. Suppose that $n\ge 1$ and put $\widehat{\widehat{X}}=X\sqcup\{y,y',z,z'\}$.   
 In the sequence    
\begin{equation}   
  \label{eq:cc12}   
\CM_x^{n-1}(\widehat{\widehat{X}})
\pile{\rTo^{c_1}\\ \rTo_{c_2}}
\CM_x^n(\widehat{X}) 
\rTo^c \CM^{n+1}_x(X)   
\end{equation}   
the compositions $cc_1$ and $c c_2$ coincide.   
   
Here  
\begin{equation} \label{eq:c1c2} 
c_1=c_{zz'} \text{ \ and \ } c_2=\theta c_{yy'}, 
\end{equation} 
 where $\theta$ is induced  by the  
involution 
of $\widehat{\widehat{X}}$  
identical on $X$ and sending $y$ to $z$ and $y'$ to $z'$.   
\end{lem}   
\begin{pf} 
Straightforward from the definition of the spaces $\CM^n_x(X)$. 
\end{pf}   
 
\medskip 
 
The following is the key technical result of this section. 
   
\subsubsection{}   
\begin{prop}{c-is-modular}   
The collection of vector spaces $\CM_x(X)$ admits a    
natural structure of a modular operad.    
\end{prop}   
   
We will construct all necessary structures by induction on degree. 
Simultaneously with checking the necessary properties we will 
establish the following characterization of  
   the kernel of the morphism~(\ref{c:def}). 
   
\subsubsection{}   
\begin{lem}{lem:ker}  
Let $k\ge 0$ and $x\in X$. The kernel of the map    
$$  
c: \CM^k_x(\widehat{X}) \to \CM^{k+1}_x(X)  
$$ 
is generated by the following three types of elements:   
   
(i)    \qquad $\alpha + \sigma_{yy'}(\alpha),$  \\ 
where 
$\alpha\in\CM^k(\widehat{X})$   and $\sigma_{yy'}$ is 
the automorphism of $\widehat{X}$  interchanging   $y$ and $y'$;   
   
(ii)   \qquad    
$\alpha\circ_{zz'}\beta -f_*(\alpha\circ_{yy'}\beta),  $ 
\\ 
where $X=A^\circ \sqcup B^\circ$;  
 \ 
$A=A^\circ \sqcup \{y,z\}$;  
 \ 
$B=B^\circ \sqcup \{y',z'\}$;  
 \  
$\alpha \in \CM^{k_1}(A)$;  
 \  
$\beta \in \CM^{k_2}(B)$;  
$k=k_1+k_2$;   
and  
$  f : X\sqcup \{z,z'\} \to  X\sqcup \{y,y'\}   
$   
is the bijection 
identical on $X$ and sending $z$ to $y$ and $z'$ to $y'$; 
   
(iii)  
\qquad 
$c_1(\alpha)-c_2(\alpha),$ \\ 
where $\alpha \in \CM_x^k(X\sqcup \{y,y',z,z'\})$, \ $k\ge 1$, and 
the maps $c_1$ and $c_2$ are given in~(\ref{eq:c1c2}). 
\end{lem}   
\medskip

It is easy to see that any element of one of the types~(i)--(iii)  
belongs to $\Ker(c)$. For type~(i) this is so because $c$ is symmetric,  
for type~(ii) it    follows from the invariance 
property~(\ref{eq:invariance}) of $c$, 
and for~(iii) this is the statement~2 of~\Lem{properties-of-c}. 
   
Therefore, to prove \Lem{lem:ker} it remains to show that any element of
$\Ker(c)$ is a combination of elements of these three  types. 
We will prove  this by induction on $k$. Before making the $k$-th step we
will identify   all the spaces $\CM^k_x(X)$ for different $x\in X$. 
\medskip

\subsubsection{Induction hypothesis}   
\label{data}   
Let us assume that   
we have the  necessary structure on $\CM^k$ for all $k \le n$.    
This includes the following components. 
\medskip

1. A canonical identification of $\CM^k_x(X)$  for different $x\in X$   
(which makes $\CM^k(X)$ well-defined) for every $k\le n$.
This means that a compatible collection of isomorphisms 
$$ \phi_*:\CM^k_x(X)\to\CM^k_{\phi(x)}(X)$$ 
is given for each automorphism $\phi:X\to X$. 
 
     We assume that  
the maps~(\ref{c:def})  \  $c:\CM^k(\widehat{X})\to\CM^{k+1}(X)$  
are equivariant with respect to automorphisms of $X$ for $k+1\le n$.   
\medskip
   
2. A collection of functorial operations   
\begin{equation} 
\circ_f: \CM^p(X) \otimes \CM^q(Y) \to \CM^{p+q+d-1}(X\sqcup Y \setminus   
(U\sqcup V)),   
\end{equation} 
defined for bijections $f:U\iso V$, where $d=|U|=|V|\ge 1$ and  
$p+q+d-1 \le n$.     
These operations satisfy conditions (i)--(iii) of
definition~\ref{ss:mod-operad} when all the compositions make sense.     
\medskip

3.  The statement of \Lem{lem:ker} is valid for all $0\le k<n$.   
\medskip

After completing the induction step we will have  
constructed this data for $k=n+1$.    
 
\medskip   
   
        Let us verify the base of     induction.   
   
\subsubsection{Tree level case:  $n=0$}      
We have $\CM_x^0(X) = \CO(X)$, therefore there is a well-defined cyclic   
structure on $\CM^0$. In particular, it gives    a canonical identification 
of spaces $\CM_x^0(X)$ for different $x\in X$.  
The structure of a cyclic operad on $\CM^0$ also gives all 
composition operations with values in $\CM^0$.  
The statement of \Lem{lem:ker} is empty for $k<0$.   
   
\subsubsection{Induction step.}   
   
Suppose that we have the structure elements (1)--(3)  
of~\ref{data} on $\CM^k$ for $k\le n$.  
Now we are going to extend this data to $\CM^{n+1}$.   
   
First in~\ref{lem:ker.step}--\ref{seq:endlemker} we will prove the 
key~\Lem{lem:ker} for $k=n$. After that we will finish the induction step 
in~\ref{acc.step}.    
   
\subsubsection{Beginning of the induction step for \Lem{lem:ker}}   
\label{lem:ker.step}   
   
In order to prove \Lem{lem:ker} for $k=n$ we first 
define a collection of functors   
$$
\CT^k:\mathbb{S}^*\to\Vect
$$ 
by the formula 
$$   
\CT^k(X)=   
\begin{cases}   
     0 & \mathrm{\ if \ } k>n+1,\\   
  \CM^k(X) & \mathrm{\ if\ } k\le n,\\   
  \CM^n(\widehat{X})/ R   & \mathrm{\ if \ } k=n+1,   
\end{cases}   
$$   
where $R$ is the subspace of $\CM(\widehat{X})$ generated   
by the elements (i) --- (iii) of~\ref{lem:ker}.

\subsubsection{} 
\begin{lem}{lem:modopt} 
The set of functors 
$\CT=\{\CT^k\}$ admits a natural structure of a modular operad.   
\end{lem} 
\begin{pf}   
We will use here the definition of modular operad 
given by \Prop{mod-operad2}. Functoriality of $\CT(X)$ with respect to $X$   
follows from  its definition. It only remains to construct    the composition
and contraction operations.     
   
By the induction hypothesis, it only remains  
to define operation and contractions with values in 
$\CT^{n+1}$.   Therefore,  we need to define  
operations of the following three kinds. 
   
\begin{enumerate}   
\item[(A)] \quad $\alpha\circ \beta$ \ for \ $\alpha \in \CT^{n+1}, \quad
\beta \in \CT^0$;
   
\item[(B)]  \quad  $\alpha\circ \beta$ \ for \ $\alpha \in \CT^k, \quad \beta    
\in \CT^{n+1-k}, \quad 1\le k\le n$;
   
\item[(C)]   
\quad $c_{xy}(\alpha)$ \ for \ $\alpha\in\CT^n(X),\quad x,y\in X,\quad X\not=\{x,y\}$.   
\end{enumerate}   
   
 {\em Definition of operations of type (A)\/.}   
 
Take any $\tilde{\alpha} \in \CM^n(\widehat{X})$ representing $\alpha$   
and define $\alpha \circ \beta$ as the image of    
$\tilde{\alpha}\circ\beta$ in $\CT^{n+1}$.   
   
This is well-defined since the subspace $R$ generated by the relations    
(i) -- (iii) is stable under multiplication by $\beta$.   
\medskip   
   
 {\em Definition of operations of type (B).\/}   
 
Take $\tilde{\alpha}\in c^{-1}(\alpha) \in \CM^{k-1}(\widehat{X})$   
and define   
$\alpha \circ \beta$ to be the image of the composition   
$\tilde{\alpha}\circ \beta$ which is defined by the induction hypothesis.    
\medskip   
   
{\em Definition of operations of type (C).\/}   
 
The map $c_{xy}:\CT^n(X)\to\CT^{n+1}(X\setminus\{x,y\})$ is obtained   
from the definition of $\CT^{n+1}(X\setminus\{x,y\})$ by identification   
(any one) of $X$ with $\widehat{X\setminus\{x,y\}}$.   
The result does not depend on the identification because  
the difference is an element of type~(i) in the subspace $R$.   
 
Let us verify that these operations define a structure of a modular 
operad on $\CT$.  This amounts to checking the following properties.
\medskip

0. $c_{xy}=c_{yx}$. 

This is one of the defining properties of $c$. 
\medskip
   
1.  {\em Independent contractions commute. }  
 
This is so because their difference in $\CM^n(\widehat{X})$ 
is an element of type (iii) in the subspace $R$ and therefore 
vanishes in $\CT^{n+1}$. 
 \medskip

2. {\em Contractions commute with compositions --- property (3)    
of~\ref{mod-operad2}.}  
 
This is true because in $\CM^n(\widehat{X})$ 
the difference of the corresponding elements 
is an element of $R$ of type (ii). 
\medskip
 
3. {\em Commutativity of the composition.} 
 
Here we need to check that  
for $\alpha \in \CT^k$, $\beta \in \CT^{n+1-k}$, $ 1\le k\le n$, 
the difference 
$ 
 \tilde{\alpha}\circ\beta-\alpha\circ\tilde{\beta} 
$ 
belongs to $R$. 
       Indeed,  
this difference can be written as   
$(c_1-c_2)(\tilde{\alpha}\circ\tilde{\beta})$  
in the notation of~(\ref{eq:c1c2}) 
and therefore is an element of $R$ of type (iii). 
   
Note that  
for $\alpha\in\CT^k,\ \beta\in\CT^l$, with $1\le k \le n-l$, 
one also has   
$$ \alpha\circ\beta=c(\tilde{\alpha}\circ\beta).$$    
   
4. {\em Associativity of the composition.}

We need to show that when $\deg\alpha+\deg\beta+\deg\gamma = n+1$ 
the composition $\alpha\circ\beta\circ\gamma$ is well-defined.  
 
Since the sum of the degrees of $\alpha$, $\beta$, and $\gamma$ 
is positive, at least one of these elements has a positive degree. 
        In the case when $\deg\alpha>0$, the element
$\alpha\circ(\beta\circ\gamma) \in \CT^{n+1}$ is the image of    
$$
\tilde{\alpha}\circ(\beta\circ\gamma)=   
(\tilde{\alpha}\circ\beta)\circ\gamma\in\CM^n.$$    
There may be two possibilities: $\deg\gamma=0$ and $\deg\gamma >0$. 
 
In the case when $\deg\gamma=0$, the composition 
$\alpha\circ\beta$ belongs to $\CT^{n+1}$ and  
is represented in $\CM^n$ by  
$\tilde{\alpha}\circ\beta$.  Therefore by  
definition of the operation of type (A) 
the composition $(\alpha\circ\beta)\circ\gamma$ is represented by     
$(\tilde{\alpha}\circ\beta)\circ\gamma\in\CM^n$.   
 
If $\deg\gamma>0$, one still has   
$$\alpha\circ\beta=c(\tilde{\alpha}\circ\beta)$$   
    and     by definition of the operation of type (B) we get  
the same result.    
 
The cases when $\deg\beta>0$ or $\deg\gamma>0$  
   are considered similarly. 
   
\medskip   
   
This concludes the construction of the modular operad structure on $\CT$.   
\end{pf}

\subsubsection{End of induction step for \Lem{lem:ker}}  
\label{seq:endlemker} 
Consider \PROP\  $\FPM(\CT)$  
generated by the modular operad $\CT$ 
   (see~\ref{seq:newprop}). 

By definition of $\CT$ there is a natural morphism $$\CO \to \CT$$ 
and therefore, by the universal property of $\CO^C$, 
we have a morphism of \PROP s   
$$   
              f:  \CO^C \to \FPM(\CT). 
$$   
       Morphism $f$ is a surjection --- it is an isomorphism in degree zero   
and the maps $\CT^k(\widehat{X})\to \CT^{k+1}(X)$ are surjective   
for $k\le n-1$ by~\Lem{properties-of-c},  and  for $k=n$ by definition of $\CT^{n+1}$.

On the other hand, we have a natural morphism of vector spaces   
$\CT^{n+1}(X)\to\CM^{n+1}_x(X)$ induced by the map $c$  
\ (\ref{c:def}). This morphism splits the morphism   
$\CO^C(X\setminus \{x\},x) \to \CP$ which proves \Lem{lem:ker} for $k=n$.  
   
 \subsubsection{End of the proof of Proposition~\ref{c-is-modular}}   
\label{acc.step}   
Using \Lem{lem:ker} for $k=n$, we can easily complete  
the induction step.    
 
First, since the kernel of $c:\CM^n_x(\widehat{X})\to\CM^{n+1}_x(X)$ is
invariant under the action of the group of automorphisms of $X$, it induces a
natural isomorphism between $\CM^{n+1}_x(X)$ for different $x\in X$.  
We also have to define new operations in $\CM$ with values in    
$\CM^{n+1}$, but this has already been done --- we defined them    
for the modular operad $\CT$ which coincides with $\CM$ in degrees    
$\le n+1$.    
   
\Prop{c-is-modular} is proved.                 \eproof 
     
\subsection{Proof of~\ref{t:left-adj} --- \ref{cor:inj}}   
\label{therest}   
   
\subsubsection{}   
\label{pf-coc=}   
Proof of~\Thm{coc=}.   
By construction of $\FPM(\FMC(\CO))$ we have a natural morphism    
$\CO \to\FPM(\FMC(\CO))$. The universal property of $\CO^C$ yields a map    
of \PROP s   
$$F:\CO^C \to \FPM(\FMC(\CO))$$   
that sends  the casimir of $\CO^C$ to the   
element $c \in \FPM(\FMC(\CO))(\emptyset,\{x,y\})$ corresponding to the    
identity element $I_{xy}$.    
   
The inverse map $ G:\FPM(\FMC(\CO)) \to \CO^C$ is constructed as follows.   
For a pair of sets $X$ and  $Y\ne \emptyset$ and an element 
$f\in \CO^C(X\sqcup Y\setminus \{y\},y)$ 
let $Y'$ be a disjoint from $Y$ copy of $Y\setminus \{y\}$ with
a fixed bijection $j: Y\setminus \{y\} \iso Y'$  and
define $G(f)$ as the composition   
$$   
X \rTo^{\id_X \sqcup c_Y} X\sqcup (Y \setminus \{y\}) \sqcup Y' 
\rTo^{f\sqcup \id_{Y'}} \{y\} \sqcup Y' \iso Y,    
$$   
   
This defines a collection of maps $\FMC(\CO)(X\sqcup Y) \to \CO^C(X,Y)$   
that automatically extends to a map of \PROP s   
$$   
G:\FPM(\FMC(\CO)) \to \CO^C   
$$   
inverting the map $F:\CO \to \FPM(\FMC(\CO))$.   
                  \eproof   
   
\subsubsection{Proof of~\Thm{t:left-adj}}   
\label{pf-left-adj}   
   
The composition $\CO \mapsto \FMC(\CO) \mapsto \FMC(\CO)^0$ is the identity.   
   
Let $\CM$ be a modular operad. We need to construct a map    
in the opposite direction   
$ \eta:\FMC(\CM^0) \to \CM$   
functorial in $\CM$.   
   
We construct $ \eta^k:\FMC(\CM^0)^k \to \CM^k$   
by induction on $k$. In the case $k=0$  these spaces are the same and the map   
$\eta^0$ is the identity.   
For $k\ge 0$ consider a  diagram   
$$   
\begin{diagram}   
\FMC(\CM^0)^k(\widehat{X}) & \rOnto^{c_{\FMC}} & \FMC(\CM^0)^{k+1}(X) \\   
  \dTo^{\eta^k} &  & \dDashto^{\eta^{k+1}} \\   
\CM^k(\widehat{X}) & \rTo^{c_{\CM}} & \CM^{k+1}(X)   
\end{diagram}   .
$$   
The map $c_{\FMC}$ is a surjection and   
\Lem{lem:ker} implies that $\eta^k(\Ker(c_{\FMC}))$ lies in    
$\Ker(c_{\CM})$.   
This allows us to extend uniquely  $\eta^k$ to $\eta^{k+1}$.   
                            \eproof 
   
\subsubsection{Proof of~\Thm{augmentation}}   
Let $\CM$ be a modular operad. Define $\CA=\Augm(\CM)\in\AugModOp$ as   
follows.   
   
Put $\CA(X)=\CM(X)$ if $X\ne\emptyset$. Furthermore, put    
$$\CA^0(\emptyset)=S^2(\CM^0(\{x\}))$$   
so that the only composition with values in $\CA^0(\emptyset)$,   
$$ \circ_{xy}:\CA^0(\{x\})\otimes\CA^0(\{y\})\to\CA^0(\emptyset),$$   
is defined.   
   
Finally, define $\CA^{n+1}(\emptyset)$ 
as the quotient $\CM^n(\{x,y\})/R$,
      where, as in~\ref{lem:ker.step}, the subspace $R$ is generated   
by the elements (i) -- (iii) of~\ref{lem:ker}.    
   
To define  a structure of an augmented modular operad on $\CA$ we will use    
the second definition~\ref{amod-op2}.   
The only missing operation is the contraction   
$$ c_{x,y}:\CA^n(\{x,y\})\to\CA^{n+1}(\emptyset)$$   
and we define it as the natural projection of $\CA^n(\{x,y\})$ to its   
quotient. The properties (i) -- (iii) of \Defn{mod-operad2}  
   can be easily verified.             \eproof 
\subsubsection{Proof of~\Thm{com=}}   
   
By      definition 
of $\FPMA(\FMCA(\CO))$ we have a natural morphism    
$\CO \to\FPMA(\FMCA(\CO))$.   
The universal property of $\CO^M$ gives a map of \PROP s   
$$F:\CO^M \to\FPMA(\FMCA(\CO))$$   
as follows. The map $F$ sends the casimir of $\CO^M$ to the   
element $c \in \FPMA(\FMCA(\CO))(\emptyset,\{x,y\})$ corresponding to the    
identity element $I_{xy}$. \\   
Similarly, $F$ sends the metric $b\in\CO^M(2,0)$   
to the element    
$$b \in\FPMA(\FMCA(\CO )) (\{x,y\},\emptyset)$$   
corresponding to the identity element $I_{xy}$.   
   
The inverse map is constructed as follows. We start with the collection   
of maps   
$$\FMC(\CO)(X\sqcup Y)\to\CO^C(X,Y)\to\CO^M(X,Y),\quad Y\ne\emptyset,$$   
defined in~\ref{pf-coc=},   
and extend it to a larger collection 
$$ 
G(X,Y):\FMCA(\CO)(X\sqcup Y) \to \CO^M(X,Y),
$$  
        where $X$ and $Y$ may be empty.
Namely,   if $Y=\emptyset,\ X\ne\emptyset$, we define 
       the map $G(X,\emptyset)$ as  the composition   
$$ 
\FMCA(\CO)(X)=\FMC(\CO)(X)\to\CO^M(\emptyset,X)   
\overset{s}{\to}\CO^M(X,\emptyset), 
$$   
        where  
the map $s$  is the composition with the  element  
        $$ b^{\otimes n}\in\CO^M(X\sqcup X',\emptyset), \ n=|X|=|X'|.      $$  
Finally,    
the map  
       $$ G(\emptyset,\emptyset):\FMCA(\CO)(\emptyset)
\to  \CO^M(\emptyset,\emptyset) 
$$ 
in positive degrees is defined as the composition of  
$G(\{x,y\},\emptyset)$ with the  contraction operation in 
$\FMCA(\CO)$ and $c\in\CO^M(\emptyset,\{x,y\})$.   
\eproof 
   
\subsubsection{Proof of~\ref{coc-vs-com} and~\ref{cor:inj}}   
\label{pf-coc-vs-com}   
   
Using the explicit definition~\ref{seq:a-newprop} of the \PROP\ generated   
by an augmented modular operad, we obtain~\Cor{coc-vs-com}.   
   
\Cor{cor:inj} is its immediate consequence.     \eproof 

\section{Enveloping algebras}   
\label{constr}

In this section we discuss different approaches to defining a universal
enveloping algebra of an algebra over an operad in a tensor category. We show
that in the case of Lie algebras  both definitions of an external universal
algebra are equivalent under mild assumptions on the tensor category. In this
case we also study the Hopf algebra structure and the notion of a center of
an enveloping algebra.  

In Section~\ref{chord} we will show that the algebra $\CA$ of chord
diagrams arising in the theory of Vassiliev knot invariants is isomorphic to
the enveloping algebra of the  universal Casimir Lie algebra  $\UL^C\in\Lie^C$.

\subsection{Definitions of enveloping algebras}
\label{def_env_alg}

The following is a standard definition of the enveloping algebra of 
an algebra over an operad in a tensor category (see,   e.g.,~\cite{hla},
Sect.~3).  Note that it defines an {\em internal\/} enveloping algebra, i.e.\
an associative algebra in the tensor category.
\medskip
 
\subsubsection{}
\begin{defn}{internal-enveloping} 
  Let  $\CO$ be an operad and $A$ be an $\CO$-algebra in a tensor category
  $\CC$.   An associative algebra $\CU$ in $\CC$ endowed with a 
collection of maps 
$$
u_n:\CO(n+1)\otimes A^{\otimes n}\to\CU,  $$  $n=1,2,\ldots$,
 is called  an      {\em internal enveloping algebra of} 
        $A$ if      $u_n$ is $\Sigma_n$-equivariant, the following diagram 
$$  \begin{diagram}
\CO(n+1)\otimes\bigotimes_{i=1}^n\CO(m_i)\otimes A^{\otimes m} 
& \rTo & \CO(m+1)\otimes A^{\otimes m}  \\
\dTo & & \dTo\\
\CO(n+1)\otimes\bigotimes_{i=1}^n\left(\CO(m_i)\otimes A^{\otimes m_i}\right) 
& & \\
\dTo & & \\
\CO(n+1)\otimes A^{\otimes n}   
& \rTo & \CU
\end{diagram}
$$
is commutative for all $n,\ m_1,\ldots, m_n$ and $m=m_1+\ldots + m_n$, 
and $\CU$ is universal with respect to these properties.   

An internal universal algebra of an $\CO$-algebra $A$, if it exists, is
unique up to an isomorphism and will be denoted by $\CU(\CO,A)$ or sometimes
simply by $\CU(A)$. 
\end{defn} 
\medskip

Internal enveloping algebras exist,  for example, when 
the tensor category $\CC$ admits colimits. In this case the enveloping 
algebra $\CU(\CO,A)$ 
can be described as the quotient of the internal tensor algebra  
\begin{equation}
\CT(\CO,A)=\bigoplus_{n\in\mathbb{N}} \CO(n+1)\otimes_{\Sigma_n} A^{\otimes n}
\end{equation}
by an ideal defined in a usual way.  
\medskip

Let $\CO$ be an operad and $\CC$ be a tensor category.
We give two different definitions of an {\em external enveloping algebra\/} of
an $\CO$-algebra $A$ in $\CC$. These algebras are ordinary associative
algebras in the category $\Vect$. In general they may not be isomorphic, but
in the important for us case of the operad $\CO=\Lie$ these definitions are
equivalent under some mild      conditions.         
\medskip 

Denote by $\Gamma$ the {\em global section\/} functor
\begin{equation}\label{eq:gamma}
\Gamma:\CC\to\Vect~, \quad \Gamma(X)=\Hom(\one,X)~.
\end{equation}

\subsubsection{}
\begin{defn}{ext_env1}
        Suppose  that the internal enveloping algebra $\CU(\CO,A)$ exists.
Then the {\em global enveloping algebra\/} $\EE(\CO,A)$ of $A$ is defined as 
\begin{equation}
\EE(\CO,A)=\Gamma(\CU(\CO,A))
\end{equation}
with operations induced by the structure of an associative algebra in $\CC$
on $\CU(\CO,A)$.   
\end{defn}
\medskip

Another way to construct an external enveloping algebra of $A$ is 
to start with the external tensor algebra of $A$ and then consider
an appropriate quotient.  
 
{\em The external tensor algebra\/} of an $\CO$-algebra $A$ is defined by the
following formula
\begin{equation} 
\label{external-tensor} 
T(\CO,A)=\bigoplus_{n\in\mathbb{N}} \CO(n+1)\otimes_{\Sigma_n} 
\Gamma(A^{\otimes n}). 
\end{equation} 
 
For every $m_1,\ldots, m_n\in\mathbb{N}$ the structure of an
$\CO$-algebra on $A$ gives a map
$$
   \CO(m_1)\otimes\ldots\otimes\CO(m_n)\otimes A^{\otimes m}\to A^{\otimes n},
$$ 
where $m=m_1+\ldots+m_n$.
Applying the functor $\Gamma$ we obtain a collection of maps 
$$\mu_{m_1,\ldots,m_n}:\CO(m_1)\otimes\ldots\otimes\CO(m_n)\otimes 
\Gamma(A^{\otimes m})\to\Gamma(A^{\otimes n})$$ 
and the following (non-commutative!) diagram   
\begin{equation}
\begin{diagram}
\CO(n+1)\otimes\bigotimes_{i=1}^n \CO(m_i) \otimes\Gamma(A^{\otimes m})
& \rTo & \CO(m+1)\otimes\Gamma(A^{\otimes m}) \\
\dTo_{id\otimes\mu_{m_1,\ldots,m_n}} & & \dInto \\
\CO(n+1)\otimes\Gamma(A^{\otimes n}) & \rInto & T(\CO,A)
\end{diagram}
\end{equation}

\medskip

\subsubsection{}
\begin{defn}{ext_env2}
The {\em external enveloping algebra\/} $U(\CO,A)$ is the quotient
of the external tensor algebra~(\ref{external-tensor}) modulo the weakest
equivalence relation making the above diagrams commutative 
for all $n$, $m_1,\ldots, m_n$.   
\end{defn}

Unlike $\EE(\CO,A)$ which is defined only when the internal enveloping
$\CU(\CO,A)$ exists, the external enveloping algebra $U(\CO,A)$ is defined
for all $\CO$, $\CC$, and $A$.

\subsubsection{}
The enveloping algebras $\EE(\CO,A)$ and  $U(\CO,A)$ are connected as follows.

Let $A$ be an $\CO$-algebra in a tensor category $\CC$. 
For $n=1,2,\ldots$ consider a map
$$ 
\CO(n+1)\otimes_{\Sigma_n}\Gamma(A^{\otimes n})\to\EE(\CO,A),
$$
 defined as the composition
\begin{eqnarray*}
\CO(n+1)\otimes_{\Sigma_n}\Gamma(A^{\otimes n})&\to&
\Gamma\left(\CO(n+1)\otimes_{\Sigma_n}A^{\otimes n}\right)\\
&\to&\Gamma(\CT(\CO,A))\to  \Gamma(\CU(\CO,A))=\EE(\CO,A) ~.
\end{eqnarray*}

These maps are compatible with the structure maps defining operations in
$\CO$ and $A$ and therefore we have a canonical homomorphism
\begin{equation} 
\label{eps} 
\epsilon_A: U(\CO,A)\to\EE(\CO,A)~.
\end{equation}

\subsubsection{} 
Both external enveloping algebras $\EE(\CO,A)$ and  $U(\CO,A)$
act functorially on all $(\CO,A)$-modules. 
 
Let $M$ be an $(\CO,A)$-module. For an element $u\in U(\CO,A)$ choose a 
representative
$\tilde{u}\in\CO(n+1)\otimes\Gamma(A^{\otimes n})$ and 
define an endomorphism $\rho$ of $M$ as the composition 
$$
M=\one\otimes M\rTo^{\tilde{u}\otimes\id_M}\CO(n+1)\otimes 
A^{\otimes n}\otimes M\to M, 
$$ 
where the second map is given
by the $(\CO,A)$-module structure on $M$. 
A straightforward check shows that $\rho$ does not depend on the choice of
the representative $\tilde{u}$ and defines a ring homomorphism
\begin{equation} 
\label{action2} 
\rho: U(\CO,A)\to\End(M). 
\end{equation} 
 
If the internal enveloping algebra $\CU(\CO,A)$ exists, one can 
similarly construct a canonical homomorphism 
\begin{equation} 
\label{action1} 
\EE(\CO,A)\to\End(M). 
\end{equation} 
In this case the homomorphism $\rho$ (\ref{action2}) coincides with the
composition of homomorphisms~(\ref{action1}) and $\epsilon_A$.

\subsubsection{}Let $i:\CO\to\CP$ be a map from an operad $\CO$ to a \PROP\  
$\CP$. The object $\mathbf{1}\in\CP$ 
has a natural $\CO$-algebra structure. 
We will denote this $\CO$-algebra by $\UA$. Its
external enveloping algebra $U(\CO,\UA)$ will be denoted by $U(\CO,\CP)$ to
stress the role of the \PROP\ $\CP$. 
Let now $\CC$ be a tensor category and let $A$ be a $\CP$-algebra in $\CC$. 
Then $A$ admits a  
natural structure of an $\CO$-algebra in $\CC$ and we have a canonical ring  
homomorphism 
\begin{equation} 
\label{univ.prop} 
U(\CO,\CP)\to U(\CO,A) 
\end{equation} 
induced by the compatible collection of vector space homomorphisms 
$$\Gamma(\mathbf{n})\to\Gamma(A^{\otimes n}).$$ 
 
\subsection{Lie algebras}    
        Until the end of this section we assume $\CO=\Lie$. 
In this case we will write $U(\fg)$ instead of $U(\Lie,\fg)$, and similarly 
for $\EE(\fg)$ and $\CU(\fg)$. 
 
We are going to apply the above constructions to the \PROP s $\Lie^C$ and
$\Lie^M$ responsible for Casimir and metric Lie algebras 
(see Section~\ref{m-c-algebras}).

\subsubsection{}
\begin{defn}{univalgs}
The algebras
$$ \UA^M=\mathbf{1}\in \Lie^M \quad \mathrm{and} \quad \UA^C=\mathbf{1}\in \Lie^C$$
are called respectively the {\em universal metric Lie algebra\/} and
 {\em universal Casimir Lie algebra\/}. 
Their external enveloping algebras will be denoted respectively by
\begin{equation}    
\label{notation}    
U^M=U(\Lie^M) \quad \mathrm{and} \quad U^C=U(\Lie^C).    
\end{equation}    
\end{defn}
\medskip
    
Corollaries~\ref{coc-vs-com} and~\ref{cor:inj} yield the     
following   result.
    
\subsubsection{}    
\begin{prop}{uc-vs-um}    
The isomorphism~(\ref{eq:com-vs-coc}) induces an isomorphism    
\begin{equation}    
\label{eq:um-vs-uc}    
\Lie^M(\mathbf{0},\mathbf{0})\otimes U^C\iso U^M.    
\end{equation}    
\end{prop}    
\eproof 
    
In particular, it implies the following 
\subsubsection{}    
\begin{cor}{cor:u-inj}    
The homomorphism    
\begin{equation}    
\label{eq:u-inj}    
i:U^C\to U^M    
\end{equation}    
induced by the functor $\Lie^C\to\Lie^M$ (\ref{eq:icm}) is injective.     
\end{cor}    
\eproof

\subsubsection{}
In order to study relations of the universal enveloping
algebras $U^M$ and $U^C$ of the universal
metric and Casimir Lie algebras $\UA^M$ and $\UA^C$
with knot invariants (see Section~\ref{chord}) we would like to have a
description of these algebras in terms of the internal enveloping algebras as
in \Defn{ext_env1}.  This seems to be impossible since the tensor categories
$\Lie^C$ and $\Lie^M$  do not admit colimits. However, this difficulty can be
resolved following  an idea of~\cite{ke2}. Since  the construction
in~\cite{ke2} misses  some important details, we present a detailed
construction of the  algebra $\CU(\UA)$ here. 
 
Let $i:\Lie\to\CP$ be a morphism from $\Lie$ to a \PROP\ $\CP$. Following 
\cite{ke2}, we construct an internal enveloping algebra $\CU(\UA)$ 
in an certain extension of the tensor category $\CP$.

Using this construction we prove~\Thm{i.e.exists+} showing 
that under mild conditions the 
internal enveloping algebra of a Lie algebra  in a tensor category
exists and the two versions of an external enveloping algebra are canonically   
isomorphic.   
 
\subsubsection{Karoubi extension}  \label{karoubi}    
We will use the following version of the Karoubi extension 
of a linear category (see~\cite{manin}, Sect.~5).
\medskip

 Let $\CC$ be a $k$-linear category with a collection $\mathcal{X}$ of pairs
 $(X,e)$, where $X\in\Ob\CC$ and $e\in\Hom(X,X)$ is an idempotent and
 for each $X\in\Ob\CC$ the pair $(X,0)$ belongs to $\mathcal{X}$.
\medskip

\begin{defn}{def:karoubi}
The {\em Karoubi extension\/}  of $\CC$ with respect to the collection 
$\mathcal{X}$ is the category $\widetilde{\CC}$ whose   
objects are pairs $(X,e)\in \mathcal{X}$ and morphisms 
from $(X,e_X)$ to $(Y,e_Y)$ are  maps $f\in\Hom_{\CC}(X,Y)$,
such that
   $$ e_Y\circ f=0\ \mathrm{and}  \  f\circ e_X=0.$$
\end{defn}

Informally speaking, the object $(X,e)$ in the category $\widetilde{\CC}$
corresponds to the kernel of $e:X\to X$.

\subsubsection{Category $\widehat{\CP}$}    

The extension $\widehat{\CP}$ of the category $\CP$ is constructed in     
two steps.    
    
First,  we construct the Karoubi extension $\widetilde{\CP}$ of the category
$\CP$  adding kernels of the idempotents in  $k\Sigma_n\subseteq\CP(n,n)$,
$n=1,2,\ldots$. The objects of the category $\widetilde{\CP}$  correspond 
to direct summands of regular representations of the symmetric groups    
$\Sigma_n$. Note that if $e_1,e_2$ are two different idempotents of
$k\Sigma_n$ corresponding to isomorphic direct summands of the regular
representation of $\Sigma_n$, then the corresponding objects of 
$\widetilde{\CP}$ are isomorphic.

After that we add to $\widetilde{\CP}$ all direct sums of  objects so that
the irreducible objects are ``of finite type''. 
Namely, let        $$\Sigma^{\widehat{}}=\coprod\Sigma_n^{\widehat{}}$$ 
be the set of irreducible representations of the symmetric groups for
$n=1,2,\ldots$.  
Then the isomorphism
classes of irreducible objects $[V]$ in $\widehat{\CP}$ are numbered 
by  $V\in\Sigma^{\widehat{}}$.    
    
Furthermore, for $X=\oplus_{i\in I}[V_i]$ and $Y=\oplus_{j\in J}[V_j]$,    
one has    
$$
\Hom_{\widehat{\CP}}(X,Y)=\prod_{i\in I}\bigoplus_{j\in J}
\Hom_{\widetilde{\CP}}([V_i],[V_j]).$$      
    
The tensor category structure on $\CP$ can be uniquely extended    
to $\widehat{\CP}$ so that the functor $\otimes$ commutes with 
coproducts. 
    
\subsubsection{}    
To construct the enveloping algebra of $\UA$ 
we use the \PBW\ theorem which allows one to realize it as a symmetric  
algebra with a deformed multiplication.

Let $\CL(V)$ be the free Lie algebra generated by $V\in\Vect$.
The symmetrization map    
$$ s:\CS(\CL(V))\to \CU(\CL(V))=\CT(V)$$    
from the symmetric algebra of $\CL(V)$ to the tensor algebra of $V$ which    
identifies with the enveloping algebra of $\CL(V)$ is bijective  by the \PBW\ 
theorem. This fact can be interpreted as an isomorphism of polynomial
functors  
\begin{equation}    
s: \CS\circ \CL\to \CT    
\label{abstract-pbw}    
\end{equation}    
which is equivalent to an 
infinite collection of identities in the 
representation rings of symmetric groups.    

Any polynomial functor on $\Vect$ defines a functor on $\widehat{\CP}$. 
Therefore, isomorphism~(\ref{abstract-pbw}) can be considered as 
an isomorphism of polynomial functors on $\widehat{\CP}$. 
This can be  viewed as the \PBW\ 
theorem for free Lie algebras in $\widehat{\CP}$.    
   
Let     $$t:\CT(\UA)\to \CS(\CL(\UA))$$
be the inverse of the isomorphism $s$.  

The Lie algebra structure on  
$\UA\in\widehat{\CP}$    
defines a map $\CL(\UA)\to\UA $ which induces a map   
$$m:\CS(\CL(\UA))\to \CS(\UA)$$
of symmetric algebras. 

\subsubsection{}    
\begin{lem}{u=}    
\

1.
The projection
\begin{equation} \label{eq:proj}    
\pi=m\circ t:\CT(\UA)\to \CS(\UA).    
\end{equation}    
is the coequalizer of the pair    
\begin{equation}    
\label{presentation-of-u}    
\CT(\CL(\UA))
\pile{\rTo^{\pi_1}\\ \rTo_{\pi_2}}
\CT(\UA),    
\end{equation}    
where
$\pi_1$ and $\pi_2$ are obtained by applying the enveloping algebra    
functor to two maps $\CL(\CL(\UA))\to \CL(\UA)$ (note that the     
enveloping algebra functor is defined on free Lie algebras).    
    
2. 
The map $\pi$ is a coalgebra morphism with respect to Hopf algebra structures
on $\CT(\UA)$ and $\CS(\UA)$ for which the sets of primitive elements
coincide with $\UA$~. 

3. 
The restriction of $\pi$ to $\CT^1(\UA)\subset \CT(\UA)$ 
is the identity map.     
\end{lem}    
\begin{pf}    \

  1. This is an old Quillen's trick --- see~\cite{q}, Appendix B.    
Compare the pair~(\ref{presentation-of-u}) with the following pair    
\begin{equation}    
\label{presentation-of-s}    
\CS(\CL(\CL(\UA))) \pile{\rTo^{\pi_1}\\ \rTo_{\pi_2}}
\CS(\CL(\UA)).    
\end{equation}    
The pair~(\ref{presentation-of-s}) is known to have $\CS(\UA)$ as a    
coequalizer since it is obtained by applying the symmetric algebra functor    
to the split sequence    
\begin{equation}    
\label{lie-bar}    
\CL^2(\UA) \pile{\rTo^{\pi_1}\\ \rTo_{\pi_2}} \CL(\UA)\lra\UA~.    
\end{equation}    
 On the other hand, the pairs~(\ref{presentation-of-s})    
and~(\ref{presentation-of-u}) are isomorphic by the \PBW\  theorem for free 
Lie algebras~(\ref{abstract-pbw}). 
    
2. It is enough to check that the map $s: \CS(\CL(\UA))\to \CT(\UA)$     
 is a morphism of coalgebras. 
This is equivalent to some identities with polynomial functors  
which can be verified on their values on $V\in\Vect$. 
They, in turn, hold because the symmetrization map 
    $$s: \CS(\CL(V))\to \CT(V)$$ 
is an isomorphism of coalgebras for $V\in\Vect$.    
    
3.     This is obvious.    
\end{pf}    
    
\subsubsection{}    
\begin{cor}{E-is-env}(\PBW\ theorem.)    
The internal enveloping algebra of $\UA\in\widehat{\CP}$ exists in
$\widehat{\CP}$. It is canonically identified with the symmetric algebra   
     $\CS(\UA)$ with multiplication induced from the multiplication in $\CT(\UA)$ 
via the projection~(\ref{eq:proj}).     
\end{cor}     

\eproof 
    
\subsubsection{}    
\begin{cor}{E-is-hopf}    
The internal enveloping algebra $\CU(\UA)$  
is a Hopf algebra in $\widehat{\CP}$.   
The composition     
$$\CS(\UA)\overset{s}{\lra}\CT(\UA)\to\CU(\UA)$$ 
defines the symmetrization map    
which is an isomorphism of coalgebras.    
\end{cor}     

\eproof 
 
\medskip

The construction of the internal enveloping algebra of $\UA\in\widehat{\CP}$ 
can be easily generalized. 
 
\subsubsection{} 
\begin{thm}{i.e.exists+}  
Let $\fg$ be a Lie algebra in a 
tensor category $\CC$. Suppose that $\CC$ admits infinite direct sums 
   and that the symmetric powers $\CS^n(\fg)$  exist in $\CC$. Then the internal
enveloping algebra $\CU(\fg)$ exists  in $\CC$. Moreover, if the functor
$\Gamma$ commutes with infinite  direct sums, the homomorphism~(\ref{eps}) 
          $$\epsilon_{\fg}: U(\fg)\to\EE(\fg)$$ 
is an isomorphism.
\end{thm} 
\begin{proof} 
Let $\CP=\FP(\Lie)$ be the \PROP~(\ref{freeprop})
generated by operad $\Lie$ and let $\UA\in\CP$ be the universal Lie algebra. 
According to \Cor{E-is-env}, its internal enveloping algebra $\CU(\UA)$
exists in category $\widehat{\CP}$.

Let $\widetilde{\CC}$ be the Karoubi extension of $\CC$ obtained by adding
kernels of all idempotents in $\CC$. Since $\CC$ admits infinite
direct sums, $\widetilde{\CC}$ admits infinite direct sums as well. 
Therefore, the functor defining the Lie algebra $\fg$ in $\CC$ extends to
a functor 
       $$\hat{\fg}:\widehat{\CP}\to\widetilde{\CC}.$$

 Since $\CU(\UA)$ is defined in 
$\widehat{\CP}$ 
by means of
split coequalizers, the image $\hat{\fg}(\CU(\UA))$ 
represents an internal enveloping algebra for $\fg$ in $\widetilde{\CC}$. 
Since $\CC$ is a full subcategory in $\widetilde{\CC}$ and $\hat{\fg}(\CU(\UA))$ 
is isomorphic to the symmetric algebra of $\fg$ in $\CC$,
this proves the existence of the enveloping algebra $\CU(\fg)$.

The algebra $\CU(\fg)$ is a split coequalizer of the pair 
\begin{equation} 
\label{presentation-of-u-gen} 
\CT(\CL(\fg)) \pile{\rTo^{\pi_1}\\ \rTo_{\pi_2}} \CT(\fg) 
\end{equation} 
constructed similarly to~(\ref{presentation-of-u}).  
Let us calculate the algebra $\EE(\fg)$. 
 
Applying the functor $\Gamma$ to the diagram~(\ref{presentation-of-u-gen}) 
we get    
\begin{equation}    
\Gamma(\CT(\CL(\fg))) \pile{\rTo^{\pi_1}\\ \rTo_{\pi_2}} \Gamma(\CT(\fg)).     
\end{equation}    
The coequalizer of this pair can be identified 
 with $U(\fg)$.    This proves the second part of the theorem. 
\end{proof} 
 
\medskip
    
The rest of this section is devoted to the description of Hopf algebra     
structures on $U^C$ and $U^M$ and to proving their commutativity.

\subsection{Hopf algebra structure}    

\

An internal enveloping algebra $\CU(\UA)$ admits a natural
Hopf algebra structure. This does not induce automatically
a Hopf algebra structure on $U(\CP)$ because the map     
$$\Hom(\mathbf{0},\mathbf{m})\otimes\Hom(\mathbf{0},\mathbf{n})    
\to\Hom(\mathbf{0},\mathbf{m}\otimes\mathbf{n})$$    
is far from being bijective.
However, in the two important cases
when $\CP=\Lie^C$ or $\CP=\Lie^M$
one  can define a Hopf algebra structure on $U(\CP)$.

We return to the use of ``the coordinate-free language'' 
of Section~\ref{sec:coord-free}
so that the arguments of our \PROP s are finite sets and not natural numbers.
    
Let $\CP=\Lie^C$ or $\Lie^M$. 
For each pair $(X,Y)$ of sets we will construct a map    
\begin{equation}    
  \label{eq:hopfing}    
  \delta_{XY}:\CP(\emptyset,X\sqcup Y)\to\CP(\emptyset,X)\otimes\CP(\emptyset,Y)    
\end{equation}    
       co-associative in a natural sense.    
    
\subsubsection{Case \ $\CP=\Lie^C$}    
\label{case-c}    
\

According to~\Thm{coc=}, we have $\CP=\FPM(\CM)$, where $\CM=\FMC(\Lie)$. Thus,
$$ 
\CP(\emptyset, X)=\bigoplus_{\substack{X=\coprod_{i\in I}X_i\\     
X_i\ne \emptyset }}    
   \bigotimes_{i\in I}\CM(X_i).$$    
    
Therefore, the left-hand side of~(\ref{eq:hopfing}) is the direct sum over 
all partitions of $X\sqcup Y$ whereas the right-hand side is the sum over 
all decomposable partitions (a union of a partition of $X$ and a partition 
of $Y$).    
The map $\delta_{XY}$ can now be defined as a natural projection.    
    
\subsubsection{Case $\CP=\Lie^M$}    \

According to~\Thm{com=}, we have $\CP=\FPMA(\CM)$, where     
$\CM=\FMCA(\Lie)$. Thus,    
$$ \CP(\emptyset, X)=S(\CM(\emptyset))\otimes    
\bigoplus_{\substack{X=\coprod_{i\in I}X_i\\     
X_i\ne \emptyset }}    
   \bigotimes_{i\in I}\CM(X_i) $$    
and the map $\delta_{XY}$ for $\CP=\Lie^M$  can be defined 
as the product of the comultiplication    
$$ S(\CM(\emptyset))\to S(\CM(\emptyset))\otimes S(\CM(\emptyset))$$    
in the symmetric algebra and the  map
constructed in~\ref{case-c}
    
\subsubsection{}    
In the cases when $\CP$ is either $\Lie^C$ or $\Lie^M$,
one finally defines the comultiplication 
\begin{equation}    
  \label{eq:comult}    
  \Hom(\mathbf{0},\CU(\UA))\to   
  \Hom(\mathbf{0},\CU(\UA))\otimes\Hom(\mathbf{0},\CU(\UA))    
\end{equation}    
on $U(\CP)$ as the composition of 
$$
  \Hom(\mathbf{0},\CU(\UA))\to\Hom(\mathbf{0},\CU(\UA)\otimes\CU(\UA))
$$ 
and 
$$
\Hom(\mathbf{0},\CU(\UA)\otimes\CU(\UA)) \to    
  \Hom(\mathbf{0},\CU(\UA))\otimes\Hom(\mathbf{0},\CU(\UA)),
$$
   where the second map is uniquely defined by~(\ref{eq:hopfing}).

\subsubsection{} In order to check that the comultiplication~(\ref{eq:comult}) 
defines a Hopf algebra structure on
$U(\CP)$, we have to verify that it is an algebra homomorphism. This is
guaranteed by  the following {\em Cake lemma\/} claiming that 
two different ways of cutting a cake $X\cup Y\cup Z\cup T$ into four pieces give 
the same result.    
    
\subsubsection{}    
\begin{lem}{cake}Let $\CP$ be $\Lie^C$ or $\Lie^M$ and let
$$\gamma_{X,Y}:\Gamma(X)\otimes\Gamma(Y)\to\Gamma(X\otimes Y)$$    
be the  natural morphism.    
The following diagram    
$$    
\begin{diagram}    
\Gamma(X\otimes Y)\otimes\Gamma(Z\otimes T) &    
  \rTo^{\delta_{XY}\otimes\delta_{ZT}}  &    
\Gamma(X)\otimes\Gamma(Y)\otimes\Gamma(Z)\otimes\Gamma(T)\\    
  \dTo^{\gamma_{X\otimes Y,Z\otimes T}} &  &    
  \dTo^{\id\otimes\sigma_{\Gamma(Y),\Gamma(Z)}\otimes\id}\\    
\Gamma(X\otimes Y\otimes Z\otimes T) &  &    
\Gamma(X)\otimes\Gamma(Z)\otimes\Gamma(Y)\otimes\Gamma(T)\\    
  \dTo^{\Gamma(\id\otimes\sigma_{YZ}\otimes\id)} &  &    
  \dTo^{\gamma_{XZ}\otimes\gamma_{YT}} \\    
\Gamma(X\otimes Z\otimes Y\otimes T) &    
  \rTo^{\delta_{X\otimes Z,Y\otimes T}}  &    
\Gamma(X\otimes Z)\otimes\Gamma(Y\otimes T)    
\end{diagram}    
$$    
is commutative for each $X,Y,Z,T\in\CP$. Here $\sigma$ denotes the    
commutativity constraint.    
\end{lem}     
       \eproof
    
\subsubsection{} It is worthwhile to give the following description    
of the sets of primitive elements in the algebras $U^C$ and $U^M$.    
    
Let $\CP$ be $\Lie^C$ or $\Lie^M$ and let $\CN$ be $\FMC(\Lie)$
       if  $\CP=\Lie^C$ and $\FMCA(\Lie)$   if $\CP=\Lie^M$. 
    
From the definition of the coproduct in $U(\CP)$ it follows
that the image of the composition 
\begin{equation}    
\label{from-m}    
\sum_{n\in\mathbb{N}}\CN(n)\to    
\sum_{n\in\mathbb{N}}\CP(0,n)=    
\sum_{n\in\mathbb{N}}\Lie(n+1)\otimes_{\Sigma_n}\CP(0,n)\to U(\CP)
\end{equation}    
consists of primitive elements. The converse is also true.    
    
\subsubsection{}    
\begin{prop}{prim=}    
The sets of primitive elements in $U^C$ and in $U^M$ coincide with the    
image of~(\ref{from-m}).    
\end{prop}    
\begin{pf}    
Enveloping algebras $U^C$ and $U^M$ are  cocommutative    
Hopf algebras. They are connected since the internal enveloping algebras    
 $\CU(\UA)$  are. They are generated as algebras by the image
of~(\ref{from-m}) which consists of primitive elements. This implies that the
image  of~(\ref{from-m}) gives {\em all} primitive elements.      
\end{pf}    
    
\subsection{Centers} \label{sec:centers}   
\subsubsection{}    
\begin{defn}{center}    
Let $A$ be an associative algebra in a tensor category $\CC$. 
An element $z\in\Gamma(A)=\Hom(\one,A)$ is called {\em central\/}    
in $A$ if the following diagram    
$$    
\begin{diagram}    
A & \rTo & \one\otimes A & \rTo^{z\otimes\id} & A\otimes A \\    
\dTo & & & & \dTo^{\text{mult.}} \\    
A\otimes\one & \rTo^{\id\otimes z} & A\otimes A & \rTo^{\text{mult.}} & A     
\end{diagram}    
$$    
is commutative.
 
The collection of all central elements of $A$ is called the {\em    
center\/} of $A$. The center is a commutative subalgebra of $\Gamma(A)$   
which we denote $Z(A)$.   
\end{defn}    
    
\subsubsection{}    
\begin{rem}{}   
Our notion of a center of an associative algebra  in a tensor category   
is one of several possible.   
   
For example, let $\CC$ be the category of super vector spaces. If $A$ is an    
associative super algebra, $\Gamma(A)$ is its zero component and the center   
we defined  is the degree zero part of (super) center of $A$. 
It coincides neither with the maximal sub-superalgebra commuting with $A$     
nor with the center of $\Gamma(A)$.    
\end{rem}

\subsubsection{}    
\begin{prop}{a-to-center}
Let $\CC$ be a tensor category admitting infinite direct sums    
and let
$\fg\in\CC$ be a Casimir Lie algebra in $\CC$, such that the 
symmetric powers $S^n(\fg)$ exist. Then the image of the map 
$$
U^C\to\Gamma(\CU(\fg))=\EE(\fg)
$$    
induced by the structure tensor functor $\fg:\Lie^C\to\CC$
belongs to the center $Z(\CU(\fg))$.    
\end{prop}    
\begin{proof}    
The category $\Mod(\fg)$ of representations of $\fg$ is a tensor     
category with direct sums and $\fg$ endowed with the adjoint action    
is a Casimir Lie algebra in $\Mod(\fg)$.     
The internal enveloping algebra of $\fg$ in $\Mod(\fg)$ is just $\CU(\fg)$    
endowed with the adjoint action of $\fg$. Then the external enveloping  
algebra  is precisely     
$$\Hom_{\Mod(\fg)}(\one,\CU(\fg))=Z(\CU(\fg)).$$    
\end{proof}

\section{Vassiliev invariants and Lie algebras}    
\label{vass}

Here we  review some facts about Vassiliev knot invariants,   
the algebra of chord diagrams, and their relationship with Lie algebra-type    
structures. For more details see~\cite{BN,Konts,Vn}.    
    
\subsection{Singular knots and  chord diagrams}   
    
A {\em singular knot\/} is an immersion $K: S^1 \rightarrow    
\mathbb{R}^3$ with     
a finite number of double self-intersections with distinct tangents.    
Framed singular knots are defined similarly. Let $\CK_n$ denote the   
set of all singular (framed) knots with $n$ double points; in particular,   
$\CK_0$ is the set of ordinary (non-singular) knots.   
   
A {\em chord diagram\/} of order $n$ is an oriented circle with $n$    
disjoint pairs of points ({\em chords\/}) on it up to an orientation    
preserving diffeomorphism of the circle.  Denote by $\CD_n$ the set of    
all chord diagrams with $n$ chords.    
    
Every singular knot $K\in \CK_n$ has a chord diagram $ch(K)\in \CD_n$    
whose chords are the inverse images of the double points of $K$.

Every (framed) knot invariant $I$ with values in an abelian    
group $k$ extends to an invariant of singular knots by the rule    
\begin{equation}    
I(K_0) = I(K_+) - I(K_-),    \label{eq:vasrel}    
\end{equation}    
where $K_0$, $K_+$, and $K_-$ are singular knots which differ only    
inside a small ball as shown on the figure below:    
    
%
%
    
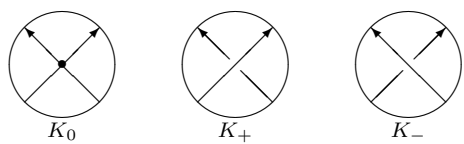
\def\KZero{    
\begin{picture}(2,2)(-1,-1)    
\put(0,0){\circle{2}}    
\put(-0.707,-0.707){\vector(1,1){1.414}}    
\put(0.707,-0.707){\vector(-1,1){1.414}}    
\put(0,0){\circle*{0.15}}    
\end{picture}}    
    
\def\KPlus{    
\begin{picture}(2,2)(-1,-1)    
\put(0,0){\circle{2}}    
\put(-0.707,-0.707){\vector(1,1){1.414}}    
\put(0.707,-0.707){\line(-1,1){0.6}}    
\put(-0.107,0.107){\vector(-1,1){0.6}}    
\end{picture}}    
    
\def\KMinus{    
\begin{picture}(2,2)(-1,-1)    
\put(0,0){\circle{2}}    
\put(-0.707,-0.707){\line(1,1){0.6}}    
\put(0.107,0.107){\vector(1,1){0.6}}    
\put(0.707,-0.707){\vector(-1,1){1.414}}    
\end{picture}}    
    
\begin{displaymath}    
\mathop{\KZero}_{K_0}\qquad    
\mathop{\KPlus}_{K_+}\qquad    
\mathop{\KMinus}_{K_-}\label{eq:VasRel}    
\end{displaymath}    

A knot invariant $I$ is called an {\em invariant of order\/} ($\le$)    
$n$ if $I(K)=0$ for any $K\in \CK_{n+1}$.    
 
We fix $k$ and denote by $\CV_n$ the set of all $k$-valued invariants of order $n$.
We have an obvious filtration    

\begin{displaymath}    
\CV_0 \subset \CV_1 \subset \CV_2 \ldots    
\subset \CV_n \subset \ldots~.    
\end{displaymath}    
    
Elements of    
$$   
 \CV = \bigcup\limits_n \CV_n   
$$   
are called {\em invariants of finite type\/} or {\em Vassiliev invariants}.     
\medskip    
    
The definition of Vassiliev invariants implies 
that the value of an invariant $I \in\CV_n$ on 
a singular knot $K$ with  $n$ self-intersections depends only on the 
diagram  $ch(K)$ of $K$.  In other words, $I$ descends to a function 
on $\CD_n$ which we still  denote by $I$.    
These functions satisfy the following relations   
\begin{equation}    
I\left(    
\Picture{\DottedCircle\FullChord[1,8]\Arc[2]\FullChord[5,9]}    
\right) -    
I\left(    
\Picture{\DottedCircle\FullChord[1,9]\Arc[2]\FullChord[5,8]}    
\right) +    
I\left(    
\Picture{\DottedCircle\FullChord[2,5]\Arc[8]\FullChord[1,9]}    
\right) -    
I\left(    
\Picture{\DottedCircle\FullChord[1,5]\Arc[8]\FullChord[2,9]}    
\right) = 0.\label{eq:4term}    
\end{equation}    
\     
    
\medskip    
    
Here each of the four diagrams has $n$ chords, but   
only chords whose endpoints lie on the solid    
arcs are shown explicitly.  The remaining $n-2$ chords 
have endpoints on the dotted arcs and are the same in all four diagrams.
    
A function $W: {\cal D}_n \rightarrow k$ is called a {\em weight    
system of order $n$\/} if it satisfies the    
{\em four-term relations\/}~(\ref{eq:4term}).    
Denote by $\CW_n$ the set of all weight systems of order $n$.    
    
Let ${\cal{A}}_n$  be the dual space to    
$\CW_n$, i.e.  the space of formal linear    
combinations of diagrams from ${\CDD_n}$ modulo the relations    
\begin{equation}    
\Picture{\DottedCircle\FullChord[1,8]\Arc[2]\FullChord[5,9]}    
 -    
\Picture{\DottedCircle\FullChord[1,9]\Arc[2]\FullChord[5,8]}    
 +    
\Picture{\DottedCircle\FullChord[2,5]\Arc[8]\FullChord[1,9]}    
 -    
\Picture{\DottedCircle\FullChord[1,5]\Arc[8]\FullChord[2,9]}    
 = 0.\label{eq:4terma}    
\end{equation}    
\\[8pt]    
    
A Vassiliev invariant  of order $n$ defines a weight system   
of order $n$ and  it is easy to see that the natural map    
 $\CV_n/\CV_{n-1} \rightarrow \CW_n$ is injective.

The remarkable fact proved by Kontsevich~\cite{Konts} is that    
if $k \supset \mathbb{Q}$ this  map is also surjective.     
In other   words, each weight system of order $n$ is a restriction to    
 ${\cal{D}}_n$ of some Vassiliev invariant.    
To prove this Kontsevich constructed a knot invariant    
$$Z:\CK_0\to \widehat{\CA} $$   
where   
\begin{equation}\label{eq:ahat}   
\widehat{\CA} = \prod_{n\ge 0}\CA_n.    
\end{equation}   
The invariant $Z$ is called {\em Kontsevich's integral}.

If $k$ is a commutative ring, then the    
product of two Vassiliev invariants of orders $m$ and $n$    
is a Vassiliev invariant of order $m+n$, therefore $\CV$    
is a filtered algebra.    
    
The space      
$$    
\bigoplus\limits_{n \ge 0}\CW_n =    
\bigoplus\limits_{n \ge 0}\CV_{n+1}/\CV_n    
$$    
becomes the associated graded algebra of $\CV$    
which induces a coproduct $\Delta$ 
on   
$$   
\CA=\bigoplus_{n\ge 0}\CA_n   
$$             defined as follows.
 
Let an element $[D]\in\CA$ be presented by a chord diagram 
$D\in\CD_m$. Then  
   
$$   
\Delta([D]) = \sum_{D_1\sqcup D_2= D } [D_1] \otimes [D_2],   $$   
where the sum is taken over all presentations of $D$ as a disjoint 
union of two subdiagrams $D_1$ and $D_2$. 
 \medskip

The operation of connected sum of diagrams induces on the coalgebra    
$\CA$ a product which makes $\CA$ a commutative and cocommutative graded
Hopf algebra.  It is called the {\em algebra of chord diagrams.}   
    
\subsection{Algebra $\CA$ and Feynman diagrams}    \label{sec:a}
\    
    
There exists an alternative description of the algebra $\CA$ of chord   
 diagrams in terms of  graphs.     
\medskip    
    
\subsubsection{}\begin{defn}    
A {\em Feynman diagram\/} of order $p$ is a graph with $2p$ vertices of    
degrees $1$ or $3$, such that each connected component has at    
least one vertex of degree $1$ and cyclic orderings are fixed on    
the set of its univalent     
({\em external\/}) vertices and on each set of three edges meeting at a     
trivalent ({\em internal\/}) vertex.\footnote{Feynman diagrams    
are called  Chinese character diagrams in~\cite{BN}, but they are indeed    
Feynman diagrams arising in the perturbative Chern-Simons-Witten    
quantum field theory.}       
Let $\CF_p$ denote the set of all    
Feynman diagrams with $2p$ vertices (up to the natural    
equivalence of graphs     
with orientations). The set $\CD_p$ of chord diagrams with $p$    
chords is a  subset of $\CF_p$.    
\end{defn}

We draw Feynman diagrams by placing their external vertices 
({\em   legs}) on a circle which is oriented  counterclockwise.    
   
  We assume that the edges meeting at each internal    
vertex are oriented counterclockwise.      
    
Denote by ${\cal{G}}_p$ the vector space generated by Feynman diagrams of order    
$p$ modulo relations    
    
\begin{equation}\label{eq:3term}    
\Picture{    
\Arc[9]    
\Arc[10]    
\Arc[11]    
\Arc[0]    
\DottedArc[1]    
\DottedArc[8]    
\thinlines    
\put(0.7,-0.7){\circle*{0.15}}    
\qbezier(-0.1,0.1)(0.2,-0.5)(0.5,0.1)    
\qbezier(0.2,-0.2)(0.3,-0.6)(0.7,-0.7)    
\put(0,-1.4){\makebox(0,0){${}_{D_Y}$}}    
}\ \ = \ \    
\Picture{    
\Arc[9]    
\Arc[10]    
\Arc[11]    
\Arc[0]    
\DottedArc[8]    
\DottedArc[1]    
\Endpoint[10]    
\Endpoint[11]    
\thinlines    
\put(0.5,-0.866){\line(-3,5){0.58}}    
\put(0.866,-0.5){\line(-3,5){0.36}}    
\put(0,-1.4){\makebox(0,0){${}_{D_{\mid\mid}}$}}    
} \ \ - \ \     
\Picture{    
\Arc[9]    
\Arc[10]    
\Arc[11]    
\Arc[0]    
\DottedArc[8]    
\DottedArc[1]    
\Endpoint[10]    
\Endpoint[11]    
\thinlines    
\put(0.5,0.1){\line(0,-1){0.97}}    
\put(0.866,-0.5){\line(-5,3){1}}    
\put(0,-1.4){\makebox(0,0){${}_{D_X}$}}    
}    
\end{equation}    
\    
    
\vspace{24pt}

More precisely, 
$$\cal{G}_p = \langle\CF_p\rangle /\langle D_Y-D_{\mid\mid}-D_X\rangle , 
$$
 where the diagrams $D_{\mid\mid}$ and    
$D_X$ are obtained from the diagram $D_Y$ by replacing its    
$Y$-fragment by the $\mid\mid$- and $X$- fragments respectively.    
    
With this notation we have the following  description of the space    
$\CA_p$ (see~{\cite{BN}}).    
    
\subsubsection{}    
\begin{prop} {prop:fd}    
    
(1) The embedding $\CDD_p \hookrightarrow \CF_p$ induces an isomorphism    
\quad $\cal{G}_p \simeq \CA_p.$    
    
(2) The following local relations hold for internal vertices in Feynman diagrams:    
    
\begin{equation}   
(i)\ \    
\begin{picture}(2,2)(0,0.375)    
\qbezier(0.5,2)(1.65,1.3)(1.5,1)    
\qbezier(1.5,2)(0.35,1.3)(0.5,1)    
\qbezier(0.5,1)(1,0)(1.5,1)    
\put(1,0){\line(0,1){0.5}}    
\end{picture}    
= -     
\begin{picture}(1.8,2)(0,0.375)    
\qbezier(0.5,2)(1,0)(1.5,2)    
\put(1,0){\line(0,1){1}}    
\end{picture}    
\ \mathrm{and}\quad    
(ii) \begin{picture}(2.5,2)(-1,-1)    
\put(0,-1.4){\line(1,1){1.4}}    
\qbezier(-0.1,0.1)(0.2,-0.5)(0.5,0.1)    
\qbezier(0.2,-0.2)(0.3,-0.6)(0.7,-0.7)    
\end{picture} =     
\begin{picture}(2.5,2)(-1,-1)    
\put(-0.034,-1.4){\line(1,1){1.4}}    
\put(0.5,-0.866){\line(-3,5){0.58}}    
\put(0.866,-0.5){\line(-3,5){0.36}}    
\end{picture} -     
\begin{picture}(2.5,2)(-1,-1)    
\put(-0.034,-1.4){\line(1,1){1.4}}    
\put(0.5,0.1){\line(0,-1){0.97}}    
\put(0.866,-0.5){\line(-5,3){1}}    
\end{picture}    
\label{eq:asjac}    
\end{equation}    
\end{prop}    
\medskip    
    
\subsection{Weight systems coming from Lie algebras}     
\label{sec:FG}

\subsubsection{}    
Here we recall a construction that    
assigns a family of weight systems to every Lie    
algebra with an invariant metric.

Let $\fg$ be a Lie algebra in a tensor category $\CC$   
with a $\fg$-invariant metric   
$b: \fg\otimes \fg \to \one$.    
To each Feynman  diagram $F$ with $m$ univalent    
vertices  we assign a tensor     
$$    
   T_{\fg}(F):\one\to \fg^{\otimes m}  
$$     
as follows.    
    
The Lie bracket $[\ ,\ ]: \fg\otimes \fg \to \fg$ can be considered    
as a tensor  
 $$\one\to\fg^*\otimes \fg^* \otimes \fg.$$
 The metric $b$ allows us to identify the $\fg$-modules     
$\fg$ and $\fg^*$, and therefore  $[\ ,\ ]$ can be considered as a tensor    
$$f: \one \to (\fg^*)^{\otimes 3}$$
and $b$ gives rise to an invariant symmetric tensor 
       $$ c: \one \to \fg\otimes \fg.$$    
    
For a Feynman  diagram $F$ denote by $T$ the set of its trivalent    
vertices, by $U$ the set of its univalent (exterior) vertices, and by    
$E$ the set of its edges. Taking $|T|$ copies of the tensor $f$ and $|E|$    
copies of the tensor $c$ we consider a new tensor    
$$    
\widetilde{T}_\fg(F) = \Bigl(\bigotimes_{v\in T} f_v\Bigr) \otimes    
                     \Bigl(\bigotimes_{\ell\in E} c_\ell\Bigr)    
$$    
which is considered as a map   
$$    
\one \to {\cal{L}}^F =     
\Bigl( \bigotimes_{v \in T}    
(\fg^*_{v,1}\otimes \fg^*_{v,2} \otimes \fg^*_{v,3})    
\Bigr)    
\otimes    
\Bigl( \bigotimes_{\ell \in E}    
(\fg_{\ell,1}\otimes \fg_{\ell,2})    
\Bigr).   
$$    
Here $(v,i), \ i=1,2,3,$ mark the three edges meeting at the    
vertex $v$ (consistently with the cyclic ordering of these edges), and    
$(\ell,j), \ j=1,2,$ denote the endpoints of the edge $\ell$.     
    
Since $c$ is symmetric and $f$ is completely antisymmetric,    
the tensor $\widetilde{T}_\fg(F)$ does not depend on the choices of    
the orderings.    
    
If $(v,i)=\ell$ and $(\ell,j)=v$, there is a natural contraction map    
         $$ \fg^*_{v,i} \otimes \fg_{\ell,j} \to \one.$$
 Composition of all such contractions gives    a map    
$$    
\gamma:    
{\cal{L}}^F \longrightarrow  \bigotimes_{u \in U} \fg = \fg^{\otimes m}, \
\text{where} \ m=|U|.     
$$    
    
The composition of $\gamma$ with $\widetilde{T}_\fg(F)$    
gives a tensor $T_{\fg}(F): \one \to \fg^{\otimes m}$.   
\medskip   
   
Often it will be convenient to    
draw Feynman diagrams with their univalent vertices along a horizontal    
line. 
 
\subsubsection{}    
\begin{exa}{}  
Let $\fg$ be a metric Lie algebra in $\Vect$ with a linear basis 
$e_1,\ e_2,\ \ldots$  
   with a metric $b$ given in this basis by $b(e_i,e_j)=b_{ij}$. 
Denote by $(b^{ij})$ the inverse of the matrix $(b_{ij})$ and 
by $f^i_{jk}$ (or $f_{ijk}$ after lowering indices by means of 
$b$) the structure constants of $\fg$ in the basis $e_1,e_2,\ldots$~. 
 
For the diagrams     
\def\DiagC{    
\begin{picture}(3,1)    
\qbezier(0.5,0)(1.5,2)(2.5,0)    
\thicklines    
\put(0,0){\line(1,0){3}}    
\end{picture}    
}    
\def\Bubble{    
\begin{picture}(4,1.5)    
\put(1.5,0){\oval(2,2)[tl]}    
\put(2.5,0){\oval(2,2)[tr]}    
\put(2,1){\circle{1}}    
\thicklines    
\put(0,0){\line(1,0){4}}    
\end{picture}}    
\def\DiagramK{    
\begin{picture}(5,1.5)    
\qbezier(0.5,0)(1.25,1.5)(2,0)    
\qbezier(3,0)(3.75,1.5)(4.5,0)    
\qbezier(1.25,0.75)(2.5,2)(3.75,0.75)    
\thicklines    
\put(0,0){\line(1,0){5}}    
\end{picture}}    
\begin{equation}    
C= \DiagC~,     
\ 
B = \Bubble~,     
\ \text{and} \    
K =\DiagramK    
\label{eq:DiagK}    
\end{equation}    
we have     
$$T_\fg(C) = \sum_{ij}  b^{ij}e_i\otimes e_j= c,    
$$     
the Casimir element corresponding to the metric~$b$, \    
$$    
T_\fg(B) = \sum b^{is}b^{tj}b^{kp}b^{lq}f_{skl}f_{pqt}e_i \otimes e_j,    
$$    
the tensor in $\fg\otimes \fg$    
corresponding to the Killing form on $\fg$ under the identification     
$\fg^* \simeq \fg$, and    
$$    
T_\fg(K) = \sum b^{in} b^{jp} b^{qr} b^{kt} b^{\ell s}    
f_{npq} f_{tsr} 
e_i \otimes e_j    
\otimes e_k \otimes e_\ell~.    
$$    
\end{exa}    
    
\subsubsection{}    The tensor $T_\fg(F)$ is invariant with respect to 
the $\fg$-action on  $\fg^{\otimes m}$ and its image $W_\fg(F)$ 
in the universal enveloping algebra $U(\fg)$ belongs to the    
center $Z(\CU(\fg))=U(\fg)^\fg$ of the {\em internal \/}    
enveloping algebra (cf.~Section~\ref{sec:centers}).   
The element  $W_\fg(F) \in Z(\CU(\fg))$ 
does not depend on the place where we cut the circle to obtain a linear
ordering of the external vertices of $F$.  
This gives a well-defined map 
$\bigoplus_{p}\langle \CF_p \rangle \to Z(\CU(\fg))$  
which vanishes on the subspace generated by the equations~(\ref{eq:asjac})
and~(\ref{eq:3term}):
relations~(\ref{eq:asjac}) follow from the anticommutativity and the Jacobi 
identity for the Lie bracket $f$, and~(\ref{eq:3term}) in this case is    
just the definition of the universal enveloping algebra as a quotient    
of the tensor algebra of $\fg$.    
    
Therefore, for every Lie algebra $\fg$ with an invariant metric   
we obtain an algebra homomorphism   
\begin{equation}    
\label{univ-wt-system}    
W_\fg: \CA \to Z(\CU(\fg))     
\end{equation}    
which is called {\em the universal weight system}    
corresponding to $\fg$.   
It is universal in the sense that the weight system $W_{\fg,R}$    
constructed using a representation $R$ of the Lie algebra $\fg$   
(see~\cite{BN})  is an evaluation of $W_\fg$:     
$$    
W_{\fg,R}(D) = \text{Tr}_{R}\bigl(W_\fg(D)\bigr).    
$$

\section{Applications}    
\label{chord}    
    
\subsection{Algebra of chord diagrams as a universal enveloping algebra}     

\
 
As an application of the results of Section~\ref{c-m-operads}    
we describe the algebra $\CA$ as the universal enveloping algebra of the   
universal Casimir Lie algebra and derive some corollaries.    
    
      \subsubsection{}    
\begin{thm}{a=}    
The algebra $\CA$ of chord diagrams     
is naturally isomorphic as  a Hopf algebra  to the external enveloping algebra  $U^C$    
of the universal Casimir Lie algebra $\UL^C$.
\end{thm}    
    
\begin{pf}    
    
Let us first construct a homomorphism $a: U^C \to \CA$. 
Define a map  $f: \Lie^C(0,n) \to \CA$ by assigning to each monomial    
in variables $c\in\Lie^C(0,2)$ and $\lambda\in\Lie^C(2,1)$     
the element of $\CA$ presented by the corresponding Feynman diagram.     
    
The map $f$ extends by linearity to the external 
tensor algebra~(\ref{external-tensor})
and gives an algebra homomorphism    
$$    
   g:T(\Lie,\Lie^C)=\bigoplus\Lie(n+1) \otimes_{\Sigma_n} 
\Lie^C(0,n) =  \bigoplus  \Lie^C(0,n) \to \CA.    
$$    
    
 One can easily  check that $g$ factors through the projection    
$$ T(\Lie,\Lie^C)\to U^C$$ 
and gives  a well-defined algebra homomorphism  $a: U^C\to\CA$.    
    
To prove that the map $a$ is also a coalgebra homomorphism we notice that     
$U^C$ is generated as algebra by its primitive elements and therefore     
it is sufficient to show that $a$ maps primitive    
elements of $U^C$ to primitive elements of $\CA$.     
According  to Proposition~\ref{prim=} the primitive elements of $U^C$ are  
images of $\Gamma(\FMC(\Lie))$. The map $a$ sends    
an element of $\Gamma(\FMC(\Lie))$ to a linear combination of connected  
graphs which are primitive in $\CA$.    
    
Since every Feynman diagram can be presented by a graph 
whose edges do not have local minima, the map $a$ is onto.    
    
It remains to show that $a$ is injective.    
Consider the diagram    
$$    
\begin{diagram}    
U^C & &\rInto^j& & U^M \\    
    &\rdOnto_a & & \ruTo_W\\    
& & \CA & &    
\end{diagram} ,    
$$    
where $W$ is the universal weight system~(\ref{univ-wt-system}) 
corresponding  
to the universal metric Lie algebra $\UA^M\in\Lie^M$, and     
$j$ is the canonical map~(\ref{eq:u-inj}).    
    
Since the diagram is commutative and $j$ is injective by~\Cor{cor:u-inj},    
this proves the injectivity of $a$.     
\end{pf}

\subsubsection{}     \label{sec:b}
Denote by $\CB$ the symmetric algebra of $\UA^C\in \Lie^C$.    
It can be described as the algebra generated by graphs analogous to
Feynman diagrams except that there is  no ordering on the set of its
univalent vertices modulo relations~(\ref{eq:asjac}).
(In the terminology of~\cite{BN} 
such graphs are called Chinese characters.)
    
\Thm{a=} and the \PBW\ theorem for Lie algebras in tensor categories   
(see \Cor{E-is-env}) give the following result of~\cite{BN}.    
\medskip
    
\begin{thm}{}    
The symmetrization map     
$$  
   \CB \rTo^\sigma \CA    
$$    
is an isomorphism of vector spaces.    
\end{thm}    
    
   \eproof

\subsection{Invariants from Casimir Lie algebras}    
    
\subsubsection{}    
\begin{prop}{}    
Let $\CC$ be a tensor category with infinite direct sums and $(\fg,t)$ be a
Casimir Lie algebra in $\CC$ such that the symmetric powers $S^n(\fg)$ exist
in $\CC$.  There is a canonical algebra homomorphism    
$$\CA\to Z(\CU(\fg)).$$    
In particular, any linear functional on the center of $\CU(\fg)$ gives    
rise to a sequence of knot invariants.    
\end{prop}    
\begin{proof}    
It follows from~\Prop{a-to-center} and \Thm{a=}.
\end{proof}   
   
One should not expect this to give new knot invariants by applying
traces of finite-dimensional representations of $\fg$ for $\CC=\Vect$.     
However, there might exist nontrivial infinite dimensional Casimir Lie    
algebras with a known center --- this may give something new.    
    
\subsubsection{Invariants from metric Lie algebras}    
\Thm{a=}  
and \Prop{uc-vs-um} imply that any finite type invariant can be obtained from
          a linear functional on the center of       
the enveloping algebra of a metric Lie algebra in a tensor $k$-linear     
category $\CC$ (one can take, for instance, $\CC=\Lie^M$ and     
    $\fg=\UA^M$). This does not contradict 
the previous claim: it is possible that some invariants coming    
from a Casimir Lie algebra in the category $\Vect$ cannot be obtained from a
metric Lie algebra in $\Vect$.    
    
\subsection{Kontsevich integral via Drinfeld's quasi-Hopf algebras}    

\

Let $(\fg,t)$ be a Casimir Lie algebra in a 
tensor category $\CC$. In~\cite{dr} Drinfeld constructed a ribbon
category $\Mod^r(\fg)[[h]]$ over $k[[h]]$,
where $\Mod^r(\fg)$ is the category of rigid  $\fg$-modules. Then a version
of Reshetikhin's construction~\cite{RT} gives a universal knot invariant lying in the     
center of the category $\Mod^r(\fg)[[h]]$.

It is tempting to try to obtain Kontsevich's integral   
$$ Z:\CK_0 \to \prod_n \CA_n h^n \subset \CA[[h]]$$   
from \Thm{a=} using Drinfeld's construction.    
We cannot apply Drinfeld's construction   
directly to $\CC=\Lie^C$ since    
this category has no rigid objects. However we can take    
$\CC=\Lie^M$ and $\fg=\UA^M$. We will obtain a knot invariant $D$    
with values in the center of $\Mod^r(\UA^M)[[h]]$.     
Let  
$$ i:U^C\to Z(\Mod^r(\UA^M)) $$
be the composition    
$$U^C\rTo^j U^M\to Z(\Mod^r(\UA^M)).$$    
    
For a knot $K$ the value $D(K)$ belongs to the image of $U^C[[h]]$ in
$Z(\Mod^r(\UA^M))[[h]].$  
If the map $i$ were injective, this would give another construction of  
Kontsevich's invariant $Z$. However, we  
do not know how  
to prove the injectivity of $i$.

\subsection{Graph complex}    

\

We showed above that  the functor $$\FMC:\CycOp\to\ModOp $$      
applied to the cyclic operad $\Lie$ naturally leads to the algebra $\CA$    
of chord diagrams. 

In another interesting situation this functor gives
Kontsevich's graph complex. We will use the variant of its definition     
given by Getzler and Kapranov~\cite{gkm}.    
    
Let $\CO$ be a cyclic Koszul operad (see~\cite{gk}, 3.2).    
The operad $\CO_{\infty}$ responsible for the homotopy $\CO$-algebras     
admits a natural cyclic structure (see~\cite{gk}, 5.4).    
For every non-empty finite set $X$ the space $\CO_\infty(X)$ is    
a complex whose only cohomology (in degree zero) is canonically    
isomorphic to $\CO(X)$.    
    
The functor $\FMC$ applied to $\CO_{\infty}$ gives 
a modular operad $\FMC(\CO_{\infty})$ in the category of    
complexes. The quasi-isomorphism  $\CO_\infty \to  \CO$    
induces a morphism $\FMC(\CO_{\infty})\to\FMC(\CO)$.     
    
The operad $\CO_{\infty}$ is a free graded cyclic operad     
   generated by $\CO^{\perp}$, where  $\CO^{\perp}$    
is the quadratic cooperad dual to $\CO$ (see~\cite{gj}).    
One has $$\CO^{\perp}=(\CO^!\{-1\})^*,$$
 where $\CO^!$ is the operad Koszul dual to $\CO$.    
    
According to~\Rem{rem:freemod},    
the modular operad $\FMC(\CO_{\infty})$ 
considered without the differential is freely generated by     
$\CO^\perp = \CO^{\!}\{-1\}^*$. Together with the differential,    
this gives the graph complex corresponding to the cyclic operad $\CO$.     
    
The graph complex $\FMC(\CO_{\infty})$ can be described in terms of     
the Feynman transform of~\cite{gkm} as follows.  The cyclic operad structure
on $\CO^!$ defines  on $\CO^!\{-1\}$ a structure of an {\em anticyclic
operad\/} (see~\cite{gk}, 2.11). 
Therefore by setting $$(\CO^!\{-1\})^n(X)=0 \  \mathrm{ for } \ n>0$$ 
we can consider $\CO^!\{-1\}$ as a twisted
 modular operad with the twist   
given by the dualizing cocycle (see~\cite{gkm}, Sect.~4).    
Then the inverse Feynman transform functor $\mathsf{F}^{-1}$ sends     
$\CO^!\{-1\}$ to the modular operad $\FMC(\CO_{\infty})$.

\subsection{BGRT conjecture}\label{sec:bgrt}    
    
Following the analogy between Lie algebras and the algebras $\CA$    
and $\CB$ of chord diagrams, Bar-Natan, Garoufalidis, Rozansky, and 
Thurston~\cite{bgrt} formulated a conjecture on the explicit form of 
the algebra isomorphism between  certain modifications $\CA'$ and 
$\CB'$ of the algebras $\CA$ and $\CB$   
    (see below).\footnote{According to Bar-Natan, 
the conjecture is now ``multiply proven'', 
see~\cite{moch} and several other txts in preparation.}    
    
As we saw in~\ref{a=}, the algebras $\CA$ and $\CB$ are isomorphic
respectively to the universal enveloping algebra 
(which is commutative in this case) and the symmetric algebra of a certain
Lie algebra in a tensor category. 
    
This allows to derive  the BGRT conjecture from our     
\Thm{a=}, \Prop{uc-vs-um}, and Kontsevich's theorem on      
Duflo-Kirillov isomorphism in arbitrary rigid tensor     
category (see~\cite{k:dq}, 8.3).    
    
Denote by $\CA'$ and $\CB'$ the spaces generated by trivalent graphs 
similar to the spaces $\CA$ and $\CB$ (see Sections~\ref{sec:a} and~\ref{sec:b})
except that connected components with no univalent vertices are allowed.    
According to \Prop{uc-vs-um} they can be described in our notation as    
$\CA'=U^M$ and $\CB'=\CS(\UA^M),$ 
and we obtain the BGRT conjecture.
\subsubsection{}    
\begin{thm}{dbgrt}    
There is a natural algebra isomorphism     
$$ \CB'\to\CA' $$    
given by the Duflo-Kirillov formula.    
\end{thm}    
         \eproof

\end{document}